\documentclass[12pt]{article}

\usepackage{graphicx}

\usepackage{amssymb}

\usepackage{amsthm}

\usepackage{epstopdf}

\usepackage{amsmath}
\usepackage{amsfonts}

\def\phi{{\varphi}}

\newcommand{\CO}[2]{ \left\langle #1 , #2 \right\rangle}

\DeclareSymbolFont{AMSb}{U}{msb}{m}{n}
\DeclareMathSymbol{\N}{\mathbin}{AMSb}{"4E}
\DeclareMathSymbol{\Z}{\mathbin}{AMSb}{"5A}
\DeclareMathSymbol{\R}{\mathbin}{AMSb}{"52}
\DeclareMathSymbol{\Q}{\mathbin}{AMSb}{"51}
\DeclareMathSymbol{\I}{\mathbin}{AMSb}{"49}
\DeclareMathSymbol{\C}{\mathbin}{AMSb}{"43}

\def\be{\begin{equation}}
\def\ee{\end{equation}}
\def\ber{\begin{eqnarray}}
\def\eer{\end{eqnarray}}

\def\beq{\begin{equation}}
\def\eeq{\end{equation}}

\def\Z{{\mathbb{Z}}}

\def\IR{{\mathbb{R}}}

\newcommand{\E}[0]{ \varepsilon}

\newcommand{ \pOm}{\partial \Omega}

\newcommand{\oo}[1]{ \overline{ #1 }  }

\begin{document}

\addtolength{\textheight}{0 cm} \addtolength{\hoffset}{0 cm}
\addtolength{\textwidth}{0 cm} \addtolength{\voffset}{0 cm}

\newenvironment{acknowledgement}{\noindent\textbf{Acknowledgement.}\em}{}

\setcounter{secnumdepth}{5}

 \newtheorem{proposition}{Proposition}[section]
\newtheorem{theorem}{Theorem}[section]
\newtheorem{lemma}[theorem]{Lemma}
\newtheorem{coro}[theorem]{Corollary}
\newtheorem{remark}[theorem]{Remark}
\newtheorem{extt}[theorem]{Example}
\newtheorem{claim}[theorem]{Claim}
\newtheorem{conj}[theorem]{Conjecture}
\newtheorem{definition}[theorem]{Definition}
\newtheorem{application}{Application}

\newtheorem*{thm*}{Theorem A}

\newtheorem{corollary}[theorem]{Corollary}

\title{A new variational principle,  convexity and supercritical Neumann problems
\footnote{Both authors  are pleased to acknowledge the support of the  National Sciences and Engineering Research Council of Canada.
}}
\author{Craig Cowan\footnote{University of Manitoba, Winnipeg Manitoba, Canada, Craig.Cowan@umanitoba.ca} \quad  Abbas Moameni \footnote{School of Mathematics and Statistics,
Carleton University,
Ottawa, Ontario, Canada,
momeni@math.carleton.ca} }

\date{}

\maketitle

\vspace{3mm}

\begin{abstract} Utilizing a new variational principle that allows dealing with problems beyond the usual locally compactness structure, we study  problems with a supercritical nonlinearity   of the type
\begin{equation}\label{eq_abstract_100}
\left\{\begin{array}{ll}
-\Delta u + u= a(x) f(u), &  \mbox{ in } \Omega, \\
u>0,    &   \mbox{ in } \Omega, \\
\frac{\partial u}{\partial \nu}= 0, &   \mbox{ on } \pOm.
\end {array}\right.
\end{equation}
To be more precise,  $\Omega$ is a bounded domain in  $\IR^N$ which satisfies certain symmetry assumptions; $ \Omega$ is a domain of `$m$ revolution' ($1\leq m<N$ and the case of   $m=1$ corresponds to radial domains)  and where $ a > 0$ satisfies  compatible symmetry assumptions along with monotonicity conditions.    We find positive nontrivial solutions of (\ref{eq_abstract_100}) in the case of suitable  supercritical nonlinearities $f$ by finding critical points of $I$ where
\[ I(u)=\int_\Omega \left\{ a(x) F^* \left( \frac{-\Delta u + u}{a(x)} \right) - a(x) F(u) \right\} dx, \]  over the closed convex cone  $K_m$  of nonnegative, symmetric and monotonic functions   in $H^1(\Omega)$ where $F'=f$ and where $ F^*$ is the Fenchel dual of $F$.
 We mention two important comments: firstly that there is a hidden symmetry in the functional $I$ due to the presence of a convex function and its Fenchel dual that makes it ideal to deal with super-critical problems lacking the necessary
  compactness requirement.  
 Secondly
the energy $I$ is not at all related to the classical Euler-Lagrange energy associated with (\ref{eq_abstract_100}).     After we have proven the existence of critical points $u$ of $I$ on $K_m$ we then unitize a new abstract variational approach (developed by one of the
present authors in  \cite{Mo,Mo2}) to show these critical points in fact satisfy
  $-\Delta  u  + u = a(x) f(u)$.
   In the particular case of $ f(u)=|u|^{p-2} u$ we show the existence of positive nontrivial solutions beyond the usual Sobolev critical exponent.

\end{abstract}


\section{Introduction}

  In this paper we consider the existence of positive solutions of the Neumann problem given by
\begin{equation}\label{eq}
\left\{\begin{array}{ll}
-\Delta u + u= a(x) f(u), &  \mbox{ in } \Omega, \\
u>0,    &   \mbox{ in } \Omega, \\
\frac{\partial u}{\partial \nu}= 0, &   \mbox{ on } \pOm,
\end {array}\right.
\end{equation}
where $\Omega$ is a bounded domain in  $\IR^N$ which satisfies  certain symmetry assumption and where $a$ is a positive sufficiently smooth function which also has some symmetry and monotonicity  properties.   When $ f$ is a subcritical nonlinearity one can utilize a standard variational approach to obtain solutions of (\ref{eq}).  With this in mind our interest is in the case of $f$ a supercritical nonlinearity;
  for example $ f(u)=|u|^{p-2} u$ where $ p> 2^*:=\frac{2N}{N-2}$.   Our approach will be to use  a new variational approach, see Theorem \ref{con2} (developed in \cite{Mo2, Mo}) over a class of functions with certain monotonicity properties,  to obtain a positive solution of (\ref{eq}).  The extra monotonicity of the functions will  give us increased ranges on the  Sobolev imbeddings and this allows one to handle  suitable supercritical nonlinearities.

\subsection{Main results and symmetry assumptions on $ \Omega$}

The domains we consider are `domains of $m$ revolution' (which we define precisely below) and of course the most basic case is a radial domain.  The next level would be what are called domains of double revolution.   Our motivation to study these special domains stems from \cite{Cabre_double} where they considered domains of double revolution  in the context of the regularity of the extremal solution associated to nonlinear eigenvalue problems of the form $ -\Delta u = \lambda f(u)$ in $ \Omega$ with $u=0$ on $ \pOm$.
 We now describe these domains.  \\

\noindent
\textbf{Domains of double revolution.}
Consider writing $ \IR^N=\IR^{m} \times \IR^{n} $ where $m+n=N$ and $m,n\geq 1.$  We define the variables $s$ and $t$ by
\[ s^2:= { x_1^2 + \cdot \cdot \cdot  + x_{m}^2}, \qquad t^2:={x_{m+1}^2 + \cdot \cdot \cdot + x_N^2 }.  \]

We say that $\Omega \subset \IR^N$ is a \emph{domain of double revolution} if it is invariant under rotations of the first $m$ variables and also under rotations of the last $n$ variables.  Equivalently, $\Omega$ is of the form $ \Omega=\{ x \in \IR^N: (s,t) \in U \}$ where $U$ is a domain in $ \IR^2$ symmetric with respect to the two coordinate axes.  In fact,
\[ U= \big \{ (s,t) \in \IR^2:  x=(x_1=s, x_2=0, ... , x_m=0,  x_{m+1}=t, ... , x_n =0 ) \in \Omega \big \},\] is the intersection of $\Omega$ with the $(x_1,x_{m+1})$ plane.  Note that $U$ is smooth if and only if $\Omega$ is smooth.  We denote $\tilde{\Omega}$ to  be the intersection of $U$ with the first quadrant of $\IR^2$.
Note that given any function $v$ defined in $\Omega$, that depends only on the radial variables $ s$ and $t$,  one has
\[ \int_\Omega v(x) dx = c(m,n) \int_{\tilde{\Omega}} v(s,t) s^{m-1} t^{n-1} ds dt,\]
where $c(m,n)$ is a positive constant depending on $n$ and $m.$ Note that strictly speaking we are abusing notation here by using the same name; and we will continuously do this in this article. Given a function $v$ defined on $ \Omega$ we will write $ v=v(s,t)$ to indicate that the function has this symmetry.   We remark that generally one requires that $ m,n \ge 2$,  but in the current work we allow the case of $ m$ or $n$ equal to $1$ as well.

\begin{extt}Let $\Omega$ be the cylinder $x^2+y^2 < 1$ with $-1 < z < 1$ in $\R^3.$ Then $\Omega$ is a domain of double revolution. In fact, by letting $s^2 ={x^2+y^2}$ and $t^2=z^2$ one has that
\[\Omega=\{(x,y,z)\in \R^3; \, \, |s|,|t|<1\}.\]
\end{extt}

\noindent
\textbf{Domains of triple revolution.} For domains of tripe revolutions and higher we adopt a more uniform notation.
 Consider writing $ \IR^N=\IR^{n_1} \times \IR^{n_2} \times \IR^{n_3}$ where $n_1+n_2+n_3=N$ and define the variables $ t_i$ via
\[ t_1^2:= { x_1^2 + \cdot \cdot \cdot  + x_{n_1}^2}, \qquad t_2^2:={x_{n_1+1}^2 + \cdot \cdot \cdot + x_{n_1 +n_2}^2 }, \qquad t_3^2:={ x_{n_1+n_2+1}^2 + \cdot \cdot \cdot + x_{N}^2 }.  \]
We say that $\Omega \subset \IR^N$ is a \emph{domain of triple revolution} if it is invariant under rotations of the first $n_1$ variables and also under rotations of the middle $n_2$ variables and the last $n_3$ variables.  Equivalently, $\Omega$ is of the form $ \Omega=\{ x \in \IR^N: (t_1,t_2,t_3) \in U \}$ where $U$ is a domain in $ \IR^3$ symmetric with respect to the three coordinate axes.  In fact,
\begin{eqnarray*} U= \Big \{ (t_1,t_2,t_3)  \in \IR^3;  x=(x_1,...,x_N) \in \Omega, \text{ where }x_1=t_1, x_{n_1+1}=t_2, x_{n_1+n_2+1}=t_3 \text{ and }\\ x_i=0 \text{ for } i\not=1, n_1+1, n_1+n_2+1 \Big \},\end{eqnarray*}
 is the intersection of $\Omega$ with the $(x_1,x_{n_1+1}, x_{n_1+n_2+1})$ plane.    We denote $\tilde{\Omega}$ to  be the intersection of $U$ with the first ``sector'' of $\IR^3$. Note that given any function $v$ defined in $\Omega$, that depends only on the radial variables $t_1,t_2,t_3$ one has
\[ \int_\Omega v(x) dx = c\int_{\tilde{\Omega}} v(t_1,t_2,t_3) t_1^{n_1-1} t_2^{n_2-1} t_3^{n_3-1} d t_1 d t_2 d t_3.\] \\
for some constant $c=c(n_1,n_2,n_3).$\\

\noindent
\textbf{Domains of $m$ revolution.}   Consider writing $ \IR^N = \IR^{n_1} \times \IR^{n_2} \times \cdot \cdot \cdot \times \IR^{n_m}$ where $ n_1 + \cdot \cdot \cdot + n_m=N$ and $n_1,...,n_m\geq 1.$
We say that $\Omega \subset \IR^N$ is a \emph{domain of $m$ revolution} if it is invariant under rotations of the first $n_1$ variables, the next $n_2$ variables, ..., and finally in the last $n_m$ variables.  We define the variables $ t_i$ via
\[ t_1^2:={ x_1^2+ \cdot \cdot \cdot  + x_{n_1}^2}, \quad t_2^2:= { x_{n_1+1}^2 + \cdot \cdot \cdot  + x_{n_1+n_2}^2},\] and similar for $ t_i$ for $ 3 \le i <m$. Finally we define
\[ t_m^2:= { \sum_{k=n_1+ n_2 + \cdot \cdot \cdot  + n_{m-1} +1}^N x_k^2}.\]
We now define
\begin{eqnarray*} U= \Big \{ t   \in \IR^m;  x=(x_1,...,x_N) \in \Omega, \text{ where }x_1=t_1,\,  x_{n_1+n_2+ \cdot \cdot \cdot  +n_{k-1} +1}= t_k \text{ for  } 2\le k\leq m,\text{ and }\\ x_i=0 \text{ for } i\not=1, n_1+1, n_1+n_2+1,..., n_1+n_2+...+n_{m-1}+1 \Big \}.\end{eqnarray*}
 We define $ \tilde{\Omega} \subset \IR^m$ to be the intersection of $U$ with the first sector of $ \IR^m$. We now define the appropriate measure
\[ d \mu_m(t) = d \mu_m^{(n_1, ... ,n_m)}(t_1, ... , t_m) = \prod_{k=1}^m t_k^{n_k-1} d t_k.\]  Given any function $v$ defined in $\Omega$, that depends only on the radial variables $t_1,t_2,..,t_m$ one has
\[ \int_\Omega v(x) dx = c(n_1,...,n_m) \int_{\tilde{\Omega}}v(t)  d \mu_m(t),\]
where $c(n_1,...,n_m)$ just depends on $n_1,...,n_m.$
Given that $\Omega \subset \R^N$  is a domain of $m$ revolution with $\sum_{i=1}^m n_i=N$, let
\[ G:= O(n_1) \times O(n_2) \times ... \times O(n_m),\]  where $O(n_i)$ is  the orthogonal group in $ \IR^{n_i}$ and
 consider
 \[ H^1_G:=\left\{ u \in H^1(\Omega): gu =u \quad \forall g \in G \right\},\]  where $gu(x):=u(g^{-1} x)$.  If $u \in H^1_G,$ then
$u$ has symmetry compatibility  with $ \Omega$, ie. $u(x)$ depends on just $ t_1,...,t_m$ and we write this as  $u(x)=u(t_1,...,t_m)$ where $(t_1,...,t_m)\in \tilde \Omega.$

\begin{remark} Now that we have clarified what we mean by a domain of $m$ revolution we can further explain our results. Indeed, if $\Omega \subset \R^N$ is a domain of $m$ revolution then one can show that the problem (\ref{eq}) admits a  positive solution $u$ of the form $u(x)=v(t_1,...,t_m)$ for some function $v: \tilde \Omega \subset \R^m \to \R.$ Moreover, by  imposing the extra condition that $\tilde \Omega$  is the unit cube in $\R^m,$ we shall be able to look for solutions with certain properties that allows us to go well beyond the Sobolev critical exponent of $\R^N.$ The next definition is  the first  step toward achieving this goal.
\end{remark}

\begin{definition} \label{def_cone}  Suppose that  $\Omega$ is a domain of $ m $ revolution in $\R^N$ with
$\tilde{\Omega} = (0,1)^m:=Q_m$ and $\sum_{i=1}^m n_i=N$. We denote by $K_m,$ the set of all  nonnegative functions $u \in  H_G^1(\Omega)$ where $u=u(t_1,...,t_m)$ is increasing with respect to each component, i.e.,
\[ K_m(n_1,..., n_m):=\left\{  u \in H_G^1(\Omega): u, u_{t_k} \ge 0  \mbox{ a.e. in $ Q_m $ for $ 1 \le k \le m$ }  \right\}.   \]
To shorten the notation, when there is no confusion, we just write $K_m$ instead of $K_m(n_1,..., n_m).$
\end{definition}
Note that $K_m$ is a closed convex cone in $ H_G^1(\Omega).$ \\

\noindent
\textbf{Assumptions on $ \Omega, f$ and $a$.} We shall assume that  $\Omega$ is a  domain of  $m$  revolution in $ \IR^N$    which further satisfies
\begin{equation} \label{assume_cube}
 \tilde{\Omega} = (0,1)^m =: Q_m.
 \end{equation}

We now consider some assumptions on the nonlinearity $f$ and $a(x)$.
\begin{itemize}
\item [$A_1$:] $f\in C^1([0, \infty))$, $f(0)=f'(0) = 0$ and $f$ is strictly increasing.

 \item [$A_2$:] There exist $p>2$  and $C>0$ such that
\[|f(t)| \leq C (1+|t|^{p-1}), \qquad \forall t \geq 0.\]

\item [$A_3$:] There exists $\mu>2$ such that for all $t \in \R$ and $ F(t):=\int_0^{|t|} f(s)\, ds$, we have
\[ |t|f(|t|) \geq \mu F(t). \]
Also there exists $l>1$ such that $2lF(t)\leq F(lt)$ for all $t\in \R.$

\item [$A_4$:] $a(x)>0$ for  a.e. $x \in \Omega$  and  $   a \in K_m \cap L^\infty(\Omega)$.

\end{itemize}  Finally, based on the number of revolutions $m$ of the domain $\Omega \subset \R^N,$ we define the $m$ dimensional critical Sobolev exponent by
$ 2^*_m:= \frac{2 m}{m-2}$ for $ m \ge 3$ and $ 2^*_1=2^*_2=\infty$. Note that $ 2^*_m:= \frac{2 m}{m-2}$ is greater than the standard Sobolev critical exponent for   $\Omega \subset \R^N.$ Indeed, if $N \geq 3$ then $2^*_N \leq  2^*_m $ as $m\leq N.$

We now state our existence theorems regarding (\ref{eq}).
\begin{theorem}  \label{exi}   Suppose that $ \Omega \subset \R^N$ is a domain of $m$ revolution with $\sum_{i=1}^m n_i=N$ and which satisfies (\ref{assume_cube}). Suppose $A_1-A_4$ hold with $ p<2_m^*$ in $A_2$. Then problem (\ref{eq}) admits at least one  positive solution $u \in K_m$.
\end{theorem}

An immediate corollary of this is the following where $\Omega$ has the desired symmetry and where $a$ does not depend fully on all $m$ variables.

\begin{corollary}\label{exx}  Suppose $ \Omega$ is a domain of $m$ revolution which satisfies (\ref{assume_cube}).  Assume that
\[ a(t_1,t_2, ..., t_m)= a(t_1, t_2, ... , t_i)\]  for some $ 1 \le i <m$.
Suppose  $A_1-A_4$ hold with $ p<2_i^*$ in $A_2$. Then problem (\ref{eq}) admits at least one  positive solution $u \in K_m$.
\end{corollary}

\begin{extt}
Consider  the  Neumann problem
\begin{equation}\label{eq100}
\left\{\begin{array}{ll}
-\Delta u + u= |x|^\alpha|u|^{p-2}u, &   x\in B_1 \\
u>0,    &   x\in B_1, \\
\frac{\partial u}{\partial \nu}= 0, &    x\in  \partial B_1,
\end {array}\right.
\end{equation}
where $B_1$ is the unit ball  centered at the origin  in $\mathbb{R}^N$, $N \geq 3.$  Note that assumptions $A_1-A_4$ in Theorem \ref{exi} hold for all $p>2$ and  $\alpha>1,$ and  therefore problem (\ref{eq100}) has a radially increasing solution $u(|x|)$.
\end{extt}

 \begin{extt} Let $\Omega \subset \R^N$ with $N\geq 3$ be a domain of double revolution with $\tilde \Omega =(0,1)^2,$  i.e.
\[\Omega=\big \{(x_1,...,x_N)\in \R^N;\,\,  x_1^2+...+x_m^2 < 1 \text{ and } x_{m+1}^2+...+x_N^2 < 1 \big\},\]
for some $1\leq m<N.$ Let $b_1, b_2: [0,1] \to (0, \infty)$ be two functions that are  twice differentiable and increasing. Consider the following problem

\begin{equation}\label{eq_abstract}
\left\{\begin{array}{ll}
-\Delta u + u= u^{p-1} b_1(\sqrt{ x_1^2 + \cdot \cdot \cdot  + x_{m}^2})b_2(\sqrt{x_{m+1}^2 + \cdot \cdot \cdot + x_N^2 }) , &  \mbox{ in } \Omega, \\
u>0,    &   \mbox{ in } \Omega, \\
\frac{\partial u}{\partial \nu}= 0, &   \mbox{ on } \pOm,
\end {array}\right.
\end{equation}
It follows from  Theorem \ref{exi}, for each   $p>2$, problem  (\ref{eq_abstract}) has a solution $u$ of the form $u(x)=v(s,t)$ for some $v:[0,1]\times [0,1] \to \R^+$  with
\[ s= \sqrt{ x_1^2 + \cdot \cdot \cdot  + x_{m}^2}, \qquad t=\sqrt{x_{m+1}^2 + \cdot \cdot \cdot + x_N^2 }. \]
Moreover, the maps $s \to v(s,t)$ and $t \to v(s,t)$ are increasing.
\end{extt}

\subsection{Outline of approach}

 Our plan is to prove existence for (\ref{eq}) by making use of  a new abstract variational principle established recently in \cite{Mo} (see also \cite{Mo1, Mo2,  Mo4}). To be more specific,
let $V$ be a reflexive Banach space, $V^*$ its topological dual  and $K$ be a closed convex subset of $V.$
Assume that  $\Phi : V \rightarrow \mathbb{R}$ is
convex,    G\^ateaux differentiable (with G\^ateaux derivative $D\Phi(u)$) and lower semi-continuous and that
 $\Lambda: Dom (\Lambda) \subset V \rightarrow V^*$ is  a linear
symmetric  operator.
Let $\Phi^*$ be the Fenchel dual of $\Phi,$ i.e.
\[\Phi^*(u^*)=\sup \{\langle u^*, u \rangle- \Phi (u); u \in V\}, \qquad u^* \in V^*,\]
where  the pairing between  $V$ and $V^*$  is denoted by
$\langle.,.\rangle$.
 Define the function $\Psi_K: V \to (-\infty, +\infty]$ by
\begin{eqnarray}
\Psi_K(u)=\left\{
  \begin{array}{ll}
      \Phi^*(\Lambda u), & u \in K, \\
    +\infty, & u \not \in K.
  \end{array}
\right.
\end{eqnarray}
Consider the functional $I_K: V \to (-\infty, +\infty]$ defined by \begin{eqnarray*}I_K(w):= \Psi_K(w)-\Phi(w).\end{eqnarray*}
A point $u \in Dom(\Psi_K)$ is said to be a critical point of $I_K$ if $D \Phi(u) \in \partial \Psi_K(u)$ or equivalently,
\[\Psi_K(v)-\Psi_K(u)\geq \langle D\Phi(u), v-u\rangle, \qquad \forall v \in V.\]
We shall now recall the following variational principle established in \cite{Mo}.
\begin{theorem}\label{con2}
Let $V$ be a reflexive Banach space and $K$ be a closed convex subset of $V.$ Let $\Phi : V \rightarrow \mathbb{R}$ be a G\^ateaux differentiable  convex and lower semi-continuous function,
and let  the linear operator
 $\Lambda: Dom (\Lambda) \subset V \rightarrow V^*$ be symmetric  and positive.  Assume that  $ u$ is a critical point of
$I_K(w)= \Psi_K(w)-\Phi(w),$
and that  there exists $v \in K$   satisfying the linear equation,
\[\Lambda v=D \Phi(u).\]
Then $u \in K$ is a solution of the equation
\begin{equation} \label{pooko} \Lambda u=D \Phi (u).\end{equation}
\end{theorem}
Before adapting this theorem to our case we make a couple of important observations.  Firstly note that $I_K$ (even if we pick $K=V$) is not the usual Euler-Lagrange energy associated with (\ref{pooko}).  The second point is that by picking $K$ appropriately one can gain compactness;  note the smaller we pick $K$ the more manageable  $I_K$ becomes which makes proving the existence of critical points of $I_K$ easier.   But this needs to be balanced with   the second part of the Theorem \ref{con2} where we need to solve the linear equation.

We now consider our case and for the purposes of clarity we consider the special case of $f(u)=|u|^{p-2} u$ and let us assume $p>2$.  Suppose that $ \Omega$ is  a domain of $m$ revolution and  then we write (\ref{eq}) in the abstract form
\begin{equation} \label{popp}
\Lambda u=D \Phi (u),
\end{equation}
 where  $\Lambda$ is the linear operator $-\Delta + 1$ and,  $\Phi$ is a suitable G\^ateaux differentiable  convex and lower semi-continuous function.  It can be easily seen that one should choose  $\Phi$ to be
  \[ \Phi(w):= \frac{1}{p}\int_\Omega a(x) |w|^p dx.\]   One can then perform  the calculations to see that $ I$  (we are omitting our choice of $K$ for now) will be
  \[ I(w)=\frac{1}{q} \int_\Omega  a(x)^{1-q} |-\Delta w + w|^q\, dx - \frac{1}{p}\int_\Omega a(x) |w|^p \, dx \] where $q:=p/(p-1)$ is the conjugate of $p$ and where we are using the $L^2$ inner product as our $V,V^*$ duality pairing (even though we have not specified $V$ yet).  Since the nonlinearity $f$ is supercritical one is unable to find critical points of $I$ on $H^1(\Omega)$ using standard variational approaches; for instance using a mountain-pass approach.   To alleviate the problems introduced by the supercritical nonlinearity we work on the cone $K_m$. Using the monotonicity of the elements of $K_m$ one obtains improved Sobolev imbeddings theorems (see Lemma \ref{imbeddings}) and this allows us to find a critical point $u$ of $I$ on $K_m$.  To  conclude that $u$ is indeed a solution of (\ref{popp}),    we then use  Corollary \ref{con3}  (an explicit version of   Theorem \ref{con2}).

\subsection{Background  when $ \Omega$ is a ball in $\R^N.$}
   We now give a background of problems related  to (\ref{eq}) in the case of supercritical nonlinearities.    In all works that we mention $\Omega$ is given by $B_1$ (the unit ball in $ \IR^N$ centered at the origin).  We mention that there are supercritical works related to (\ref{eq}) in the case of nonradial domains but they are generally problems which contain a small parameter $ \E$.

 In  \cite{first_rad_neum} they considered the variant of (\ref{eq}) given by $ -\Delta u + u = |x|^\alpha u^p$ in $B_1$ with $ \frac{ \partial u}{ \partial \nu}=0$ on $ \partial B_1$.
 They prove the existence of a positive radial solutions of this equation with arbitrary growth using a shooting argument.   The solution turns out to be an increasing function.
They also perform numerical computations to see the existence of positive oscillating solutions.  In  \cite{serra_tilli} they considered  (\ref{eq})   along with the classical energy associated with the equation given by
  \[E(u):=\int_{B_1} \frac{ | \nabla u|^2 +u^2}{2}\, dx  -\int_{B_1} a(|x|)F(u)\,  dx,\] where $F'(u)=f(u)$.    Their goal was  to find critical points of $E$ over $H^1_{rad}(B_1):=\{ u \in H^1(B_1):  u \mbox{ is radial} \}$.  Of course since $f$ is supercritical the standard approach of finding critical points will present difficulties and hence their idea was to find critical points of $E$ over the cone $ \{ u \in H_{rad}^1(B_1): 0 \le u, \mbox{ $u$ increasing} \}$.   Doing this is somewhat standard but now the issue is the critical points don't necessarily correspond to critical points over $H_{rad}^1(B_1)$ and hence one can't conclude the critical points solve the equation.   The majority of their work is to show that in fact the critical points of $E$ on the cone are really critical points over the full space.
 In \cite{first_grossi},
  \begin{equation}\label{eq_grossi_1}
\left\{\begin{array}{ll}
-\Delta u + V(|x|) u= |u|^{p-2}u, &   \mbox{ in }  B_1 \\
u>0,    &   \mbox{ in }  B_1, \\
\end {array}\right.
\end{equation} was examined under both homogeneous Neumann and Dirichlet boundary conditions.  We will restrict our attention to their results regarding the Neumann boundary conditions.  Consider $G(r,s)$ the Green function of the operator
\[ \mathcal{L}(u)=- u'' - \frac{N-1}{r} u' + V(r) u, \quad u'(0)=0,\] with $u'(1)=0$. Define now $H(r):= (G(r,r))^{-1} | \partial B_1| r^{N-1}$ for $ 0<r \le 1$.  One of their results states that for $ V \ge 0$ (not identically zero) if $H$ has a local minimum at  $\overline{r} \in (0,1]$ then for $p$ large enough, (\ref{eq_grossi_1}) has a solution with Neumann boundary conditions and the solutions have a prescribed asymptotic behavior as $ p \rightarrow \infty$.   Additionally they can find as many solutions as $H$ has local minimums.    This work contains many results and we will list one of which more related.  For $ V=\lambda>0,$ the problem  (\ref{eq_grossi_1}) has a positive nonconstant solution with Neumann boundary conditions provided $p$ is large enough.   This methods used in  \cite{first_grossi} appear to be very different from the methods used in the all the other works.
 It appears the works of \cite{serra_tilli} and \cite{first_grossi} were done completely  independent of each other.
 The next work related to (\ref{eq}) was \cite{Weth} where they considered
 \begin{equation}\label{eq_weth}
\left\{\begin{array}{ll}
-\Delta u +b(|x|)x \cdot \nabla u+ u= a(|x|)f(u), &  \mbox{ in }  B_1 \\
u>0,    &   \mbox{ in } B_1, \\
\frac{\partial u}{\partial \nu}= 0, &    \mbox{ on }  \partial B_1,
\end {array}\right.
\end{equation} where $f$ is a supercritical nonlinearity and where various assumptions were imposed on $b$.   Their approach was similar to \cite{serra_tilli} in the sense that they also worked on the cone  $\{ u \in H_{rad}^1(B_1): 0 \le u, \mbox{ $u$ increasing} \}$ but instead of using a variational approach they used a topological approach.   They were able to weaken the assumptions needed on $f$.   In the case of $a=1$ one sees that the constant $u_0$ is a solution provided $ f(u_0)=u_0$.   In \cite{Weth} they have showed that (\ref{eq_weth}) has a positive nonconstant solution in the case of $b=0$ provided there is some $u_0>0$ with $ f(u_0)=u_0$ and $ f'(u_0)> \lambda_2^{rad}$ which is the second radial eigenvalue of $ -\Delta +I$ in the unit ball with Neumann boundary conditions.  Note that this result shows there is a positive nonconstant solution of (\ref{eq}) provided $ p-1> \lambda_2^{rad}$.
In \cite{system} they considered various elliptic systems of the form
\begin{equation*}   \label{sys_theirs}
\left\{\begin{array}{ll}
-\Delta u + u=  f(|x|,u,v), &   \mbox{ in }  B_1 \\
-\Delta v + v=  g(|x|,u,v), &   \mbox{ in }  B_1 \\
\frac{\partial u}{\partial \nu}=\frac{\partial u}{\partial \nu}= 0, &    \mbox{ on }   \partial B_1.
\end {array}\right.
\end{equation*}  In particular they examined the gradient system when $ f(|x|,u,v)=G_u(|x|,u,v),  g(|x|,u,v)=G_v(|x|,u,v)$ and they  also considered the Hamiltonian system version where $ f(|x|,u,v)=H_v(|x|,u,v),  g(|x|,u,v)=H_u(|x|,u,v)$.  In both cases there obtain positive solutions under various assumptions (which allowed supercritical nonlinearities).  They also obtain positive nonconstant solutions in the case of $f(|x|,u,v)=f(u,v)$, $ g(|x|,u,v)=g(u,v)$;  note in this case there is the added difficulty of avoiding the possible constant solutions.

These results were extended to $p$-Laplace versions in \cite{secchi}.    The methods of \cite{first_grossi} were extended to prove results regarding multi-layer radials solutions in \cite{grossi_new}.   We also mention the work of  \cite{den_serra} where problems on the annulus were considered.    We also mention the very recent works which extend some results and answer some open questions; see \cite{add1,add2,add3,add4,add5,add6}.

 In \cite{ACL} we considered
(\ref{eq}) in the case of $ f(u)=|u|^{p-1} u$.  Using a new variational principle we obtained positive solutions of (\ref{eq}); assuming the same assumptions as the earlier works.   In the case of $a(x)=1$ we obtain the existence of a positive nonconstant solution of (\ref{eq}).  We remark our approach allowed us to deal directly with the supercritical nonlinearity without the need to cut the nonlinearity off.

We mention is that there is another type of supercritical problem that  one can examine on $B_1$.   One can examine supercritical equations like (\ref{eq}) or the case of zero Dirichlet boundary conditions when $a$ is radial and $ a=0$ at the origin; a well known case of this is the H\'enon equation given by $ -\Delta u = |x|^\alpha u^p $ in $B_1$ with $u=0$ on $ \partial B_1$ where $0<\alpha$.   In \cite{Ni} it was shown the H\'enon equation has a positive solution if and only if $ p< \frac{N+2 +2 \alpha}{N-2}$,  and note this includes a range of supercritical $p$.  This increased range of $p$ is coming from the fact that $a=0$ at the origin.   We mention this phenomena is very different than what is going on in the above works.   Results regarding positive solutions of supercritical  H\'enon equations on general domains have also been obtained,  see \cite{cowan} and \cite{glad}.

One final point we mention is that there has been extensive study of subcritical, critical and supercritical Neumann problems on general domains in the case of (\ref{eq}) when $a=1$ and where the equation involves a parameter that is sent to either zero or infinity.   We have not attempted to discuss this problem but the interested reader should consult, for instance, \cite{100,50,101,102,103,104,105,106}.

\section{Elliptic problems on domains of $m$ revolution} \label{elliptic_pr}

In this section we discuss the issue of solving equations on   domains of $m$ revolution in $\R^N$.  We begin with the standard definition of a weak solution to a Neumann boundary value problem.

\begin{definition}  We say $v$ is  a weak solution of
 \begin{equation} \label{weak_eq}
 \left\{
\begin{array}{rrl}
-\Delta v +v &=&h(x) \qquad \hfill \mbox{ in } \Omega, \\
\partial_\nu v&=& 0 \qquad \hfill  \mbox{ on } \pOm,
\end{array}
\right.
\end{equation} provided $ v \in H^1(\Omega)$ and  satisfies
\[ \int_\Omega \nabla v \cdot \nabla \eta +  v \eta \; dx = \int_\Omega h(x) \eta, \qquad \forall \eta \in H^1(\Omega).\]
\end{definition}

Given $\Omega \subset \R^N$ which is a domain of $m$ revolution with $\sum_{i=1}^m n_i=N$ and, a function $h: \Omega \to \R$ that  has symmetry compatible with $ \Omega$, i.e. $h(x)$ depends on just $ t_1,...,t_m$ (we write this as $h=h(t_1,...,t_m)$), it is natural to look for a solution of (\ref{weak_eq}) satisfying the same symmetry properties.  Recall that
\[ G= O(n_1) \times O(n_2) \times ... \times O(n_m),\]   and
 \[ H^1_G=\left\{ u \in H^1(\Omega): gu =u \quad \forall g \in G \right\},\]  where $gu(x):=u(g^{-1} x)$.     To find a solution of (\ref{weak_eq}) it is sufficient (using the principle of symmetric
criticality)  to find a critical point of
\[ E_\Omega(v):=\frac{1}{2} \int_\Omega | \nabla v|^2 + v^2 dx - \int_\Omega hv dx, \] over $H^1_G(\Omega)$;  i.e.  to find a $ v \in H^1_G(\Omega)$ such that
 \begin{equation} \label{weak_G}
 \int_\Omega \nabla v \cdot \nabla \eta + v \eta \; dx = \int_\Omega h(x) \eta \;  dx, \; \;  \forall \eta \in H_G^1(\Omega).
 \end{equation}
Note that we can identify $H^1_G(\Omega)$ with $Y_m$ where $Y_m:=\{ v: \tilde{\Omega} \rightarrow \IR : \|v\|_{Y_m} <\infty\}$ with
  \[ \|v\|_{Y_m}^2= \int_{\tilde{\Omega}} \left( \sum_{k=1}^m v_{t_k}^2 + v^2\right) d \mu_m(t), \qquad d \mu_m(t)= \prod_{k=1}^m t_k^{n_k-1} d t_k. \]
   Note we are using here that
  $ v_{x_i}= v_{t_1} \frac{x_i}{t_1} $ for $ 1 \le i \le n_1$;  $ v_{x_i} = v_{t_2} \frac{x_i}{t_2}$ for $n_1+1 \le i \le n_2$ and we can carry on like this.  So from this we see that $ | \nabla_x v|^2 = \sum_{k=1}^m v_{t_k}^2$.\\
     Note that $Y_m$ is not $H^1(\tilde{\Omega})$ after noting the degenerate weights in $ d \mu_m$.  Also  note that if $ v \in H^1_G(\Omega)$ satisfies (\ref{weak_G}) then given $ \eta \in H^1_G(\Omega)$ we have
 \begin{eqnarray} \label{comp}
 c(n_1,...,n_m)\int_{\tilde{\Omega}} h \eta \; d \mu_m(t) &=& \int_\Omega h \eta dx \nonumber \\
 &=& \int_\Omega \nabla v \cdot \nabla \eta + v \eta \; dx  \nonumber \\
 &=& c(n_1,...,n_m) \int_{\tilde{\Omega}}  \left( \sum_{k=1}^m v_{t_k} \eta_{t_k} + v \eta \right) d \mu_m(t),
 \end{eqnarray} for all $ \eta \in H^1_G(\Omega)$; note we are identifying $H^1_G(\Omega)$ and $Y_m$ without changing notation.   So we see that $v \in Y_m$ satisfies
\begin{equation} \label{WEAK_Y_first}
 \int_{\tilde{\Omega}}  \left( \sum_{k=1}^m v_{t_k} \eta_{t_k} + v \eta\right) d \mu_m(t) =\int_{\tilde{\Omega}} h \eta \; d \mu_m(t), \quad \forall \eta \in Y_m.
\end{equation}
\noindent
\textbf{Notation.}  For notational  convenience we set
\begin{equation} \label{notation_100}
(\nabla_t v)_k = v_{t_k} \mbox{ for $ 1 \le k \le m$}, \qquad \Delta_t v = \sum_{k=1}^m v_{t_k t_k}.
\end{equation}
 Integrating the (\ref{WEAK_Y_first})   by parts formally one sees that $v_m$ satisfies
\begin{eqnarray}
0 & = & \int_{\tilde{\Omega}} \left( h + \sum_{k=1}^m \left\{v_{t_k t_k} + \frac{n_k-1}{t_k}v_{t_k} \right\} - v \right) \eta d \mu_m(t) \\
&& - \int_{\partial \tilde{\Omega}} \eta \left( \sum_{k=1}^m v_{t_k} \nu^k \right) \prod_{i=1}^m t_i^{n_i-1},
\end{eqnarray}  where $ \nu=(\nu^1,..., \nu^m)$ is the outward pointing normal on $ \partial \tilde{\Omega}$.
From this we see that $v$ should satisfy

 \begin{equation}  \label{weak_obs}
 \left\{
\begin{array}{rrl}
   -\Delta_t v  - \sum_{k=1}^m \frac{n_k-1}{t_k}v_{t_k}  + v    &=& h(t)  \qquad \hfill \mbox{ in } \tilde{\Omega}, \\
\partial_\nu v&=& 0 \qquad \hfill  \mbox{ on } \partial \tilde{\Omega}.
\end{array}
\right.
\end{equation}    Note that  the boundary condition on the `inner boundaries' $ t_k=0$ for $ 1 \le k \le m$ is not coming  from the weak formulation of the problem but rather from the symmetry of $v$ and this requires $v$ to have sufficient regularity.\\

We now prove a result regarding monotonicity  of solutions of (\ref{weak_eq}) provided $h$ is monotonic.  This result is crucial when applying Theorem  \ref{con2}, in particular when showing if $u \in K$ is a critical point of $I_K$ over $K$  then there is some $v \in K$ which satisfies $\Lambda v=D \Phi(u).$

\begin{proposition}  \label{mono_linear}  Assume $ \Omega$ is domain of $m$ revolution in $ \IR^N$ which satisfies (\ref{assume_cube}) and suppose $ 0 \le h \in K_m \cap C^{1,1}(\oo{\Omega})$.  Then there exists a unique  $ v  \in K_m  \cap C^{2,\alpha}(\oo{\Omega})$ (any $ 0<\alpha<1$)  which satisfies (\ref{weak_eq}).
\end{proposition}
We shall  provide a proof for Proposition \ref{mono_linear} in the Appendix. We conclude this section by proving some improved imbeddings for functions in $K_m$ given in Definition \ref{def_cone}. Indeed,
working on $K_m$, will allow us to obtain improved Sobolev critical exponents that are essential   to consider supercritical problems from a variational point of view.

\begin{lemma} \label{imbeddings} (Improved Sobolev imbeddings on $K_m$) Suppose $\Omega$ is a  domain of $m$ revolution in $\R^N$  which satisfies (\ref{assume_cube}).
\begin{enumerate} \item For $ m \ge 3$ and for $ 1 \le q \le 2_m^*$ (in the case of $ m=1,2$ for $ 1 \le q < \infty$) there is some $ C_q$ such that
\begin{equation} \label{imbed_1} \| u \|_{L^q(Q_m)} \le C_q \| u \|_{H^1(\Omega)} \qquad \forall u \in K_m.
\end{equation}

\item For $ m \ge 3$ and for $ 1 \le q \le 2_m^*$ (in the case of $ m=2$ for $ 1 \le q < \infty$) there is some $ C_q$ such that
\begin{equation} \label{imbed_2} \| u \|_{L^q(\Omega)} \le C_q \| u \|_{H^1(\Omega)} \qquad \forall u \in K_m.
\end{equation}

\end{enumerate}
\end{lemma}

\noindent
\textbf{Proof.}
We first prove 1). The second part  is a direct consequence of part 1).
\begin{enumerate}  \item  Consider $A_1:= (2^{-1},1)^m$ and note we can decompose $Q_m$ into the union $2^m$ disjoint cubes which are translations of $A_1$ (say $ A_1,A_2,..., A_{2^m}$) where we are missing a set of measure zero of $Q_m$.  Let $ u \in K_m$ and let $1 \le q <\infty$ if $m=2$ and for $ m \ge 3$ let $ 1 \le q \le 2_m^*$. Then by the  $m$ dimensional Sobolev imbedding we have
\[ \| u \|_{L^q(A_1)} \le C_q \| u \|_{H^1(A_1)}.\]  Note there is some constant $ C=C(n_1,n_2,...,n_m)$ such that
$ \|u \|_{H^1(A_1)} \le C \| u \|_{Y_m}=\tilde{C} \|u\|_{H^1(\Omega)}$.   By monotonicity of $u$  we have $ \| u\|_{L^q(A_k)}^q \le \| u \|_{L^q(A_1)}^q \le \tilde{C}^q \|u \|_{H^1(\Omega)}^q$ for $ 1 \le k \le m$.  Summing in $k$ gives (\ref{imbed_1}).

\item Note that for any $ 0 \le u$ with $ u=u(t_1, ... , t_m)$ we have
\[ \|u \|_{L^q(\Omega)}^q = C \int_{Q_m} u^q d \mu_m(t) \le C \int_{Q_m} u^q d t_1,..., dt_m,\] and we can then use the part 1 of this lemma to obtain the desired result.
\end{enumerate}
\hfill $\square$\\

\section{Preliminaries from convex analysis} \label{convex}

In this section we recall some important definitions and results from   convex analysis and  minimax principles for lower semi-continuous functions.

Let $V$ be a  real Banach  space and $V^*$ its topological dual  and let $\langle .,. \rangle $ be the pairing between $V$ and $V^*.$
The weak topology on $V$ induced by $\langle .,. \rangle $ is denoted by $\sigma(V,V^*).$  A function $\Psi : V \rightarrow \mathbb{R}$ is said to be weakly lower semi-continuous if
\[\Psi(u) \leq \liminf_{n\rightarrow \infty} \Psi(u_n),\]
for each $u \in V$ and any sequence ${u_n} $ approaching $u$ in the weak topology $\sigma(V,V^*).$
Let $\Psi : V \rightarrow \mathbb{R}\cup \{\infty\}$ be a proper convex  function. The subdifferential $\partial \Psi $ of $\Psi$
is defined  to be the following set-valued operator: if $u \in Dom (\Psi)=\{v \in V; \, \Psi(v)< \infty\},$ set
\[\partial \Psi (u)=\{u^* \in V^*; \langle u^*, v-u \rangle + \Psi(u) \leq \Psi(v) \text{  for all  } v \in V\}\]
and if $u \not \in Dom (\Psi),$ set $\partial \Psi (u)=\varnothing.$ If $\Psi$ is G\^ateaux differentiable at $u,$ denote by $D \Psi(u)$ the derivative of $\Psi$ at $u.$ In this case  $\partial \Psi (u)=\{ D  \Psi(u)\}.$\\
 The  Fenchel  dual of an arbitrary function $\Psi$ is denoted by  $\Psi^*,$ that is function on $V^*$ and is defined by
\[\Psi^*(u^*)=\sup \{\langle u^*, u \rangle- \Psi (u); u \in V\}.\]
Clearly $\Psi^*: V^* \rightarrow \mathbb{R}\cup \{\infty\}$  is convex and  weakly lower semi-continuous.  The following standard  result  is crucial in the subsequent analysis (see \cite{Ek-Te, Ek2} for a proof).
\begin{proposition}\label{var-pro} Let $\Psi : V \rightarrow \mathbb{R}\cup \{\infty\}$ be an arbitrary function. The following statements hold:\\
(1) $\Psi^{**}(u) \leq \Psi(u)$ for all $u \in V.$ \\
(2) $\Psi(u)+\Psi^*(u^*)\geq \langle u^*, u \rangle$ for all $u \in V$ and $u^*\in V^*.$\\
(3) If $\Psi$ is convex and lower-semi continuous then $\Psi^{**}=\Psi$ and the following assertions are equivalent:
\begin{itemize}
\item  $\Psi (u) +\Psi^*(u^*) = \langle u, u^* \rangle.$
\item  $u^* \in \partial \Psi (u).$
\item  $ u \in \partial \Psi^* (u^*).$
\end{itemize}
\end{proposition}
The above Proposition shows that for a convex lower semi-continuous function  $\Psi$ one has \[\partial \Psi^* =\big (\partial \Psi\big)^{-1}.\]
We shall now recall some notations and results for the minimax principles of  lower semi-continuous functions.
\begin{definition}\label{crit100}
Let $V$ be a real Banach space,  $\phi\in C^1(V,\mathbb{R})$ and $\psi: V\rightarrow (-\infty, +\infty]$ be proper (i.e. $Dom(\psi)\neq \emptyset$), convex and lower semi-continuous.
A point $u\in  V$ is said to be a critical point of \begin{equation} \label{form}I:=  \psi-\phi \end{equation} if $u\in Dom(\psi)$  and if it satisfies
the inequality
\begin{equation}
 \CO{D \phi(u)}{ u-v} + \psi(v)- \psi(u) \geq 0, \qquad \forall v\in V.
\end{equation}
\end{definition}

\begin{definition}
We say that the functional $I=\psi-\phi,$ given in (\ref{form}), satisfies the Palais-Smale compactness  condition (PS)   if
every sequence $\{u_n\}$ such that
\begin{itemize}\label{2}
\item  $I[u_n]\rightarrow c\in \mathbb{R},$
\item  $\CO{D \phi(u_n)}{ u_n-v} + \psi(v)- \psi(u_n) \geq -\epsilon_n\|v- u_n\|, \qquad \forall v\in V.$
\end{itemize}
where $\epsilon_n \rightarrow 0$, then $\{u_n\}$ possesses a convergent subsequence.
\end{definition}

The following is  proved in \cite{szulkin}.
\begin{theorem}\label {MPT}
(Mountain Pass Theorem).  Suppose that
$I : V \rightarrow (-\infty, +\infty ]$ is of the form (\ref{form}) and satisfies the Palais-Smale   condition and  the Mountain Pass Geometry (MPG):
\begin{enumerate}
\item $I(0)= 0$.
\item  there exists $e\in V$ such that $I(e)\leq 0$.
\item there exist $\alpha>0$ and  $0<\rho<\|e\|$ such that  for every $u\in V$ with $\|u\|= \rho$ one has $I(u)\geq \alpha$.
\end{enumerate}
Then $I$ has a critical value $c\geq \alpha$ which is  characterized by
$$c= \inf_{g\in \Gamma}  \sup_{t\in [0,1]} I[g(t)],$$
where   $\Gamma= \{g\in C([0,1],V): g(0)=0, g(1)= e\}.$
\end{theorem}

\section{Existence results}

In this section we assume that $\Omega \subset \R^N$ is a domain of $m$ revolution with $\sum_{i=1}^mn_i=N$ and $\tilde \Omega=Q_m.$ Throughout this section we always assume that assumptions $A_1-A_4$ hold. For a  function $f:[0,\infty) \to \R$ satisfying $A_1-A_3,$ define $F : \R \to \R$ and  by
\[F(t)=\int_0^{|t|} f(s)\, ds,\]
and let $F^*: \R \to (-\infty, +\infty]$ be the Fenchel dual of $F,$ i.e.,
\[F^*(s)=\sup_{t \in \R} \{ts-F(t)\}.\]
To adapt Theorem \ref{con2} in  our case,  consider
 the Banach space $V= H^1(\Omega)\cap L^{p}(\Omega)$, where $2<p<2_m^*$
and $V$ is equipped  with the following norm
\[\|v\|:= \|v\|_{H^1(\Omega)}+ \|v\|_{L^{p}(\Omega)}.\] \\
Assuming $V^*$ is  the topological dual of $V,$ the pairing $\langle .,. \rangle $ between $V$
and $V^*$ is defined by
\[\langle v,v^* \rangle=\int_\Omega v(x) v^*(x)\, dx,\qquad \forall v \in V,\, \forall v^* \in V^*.\]
For $v\in V$ define the operator $A: Dom(A)\subset V \to V^*$ by
$Av:=  -\Delta v+v,$ where \[Dom(A)=\{v \in V; \frac{\partial v}{\partial n}=0, \quad \& \quad Av \in V^*  \}.\]
Note that one can rewrite the problem (\ref{eq}) as
$$A u=D \varphi(u),$$
where $\phi: L^p(\Omega) \to \R$ is defined by
\begin{equation} \label{gnov} \varphi(u)= \int_{\Omega} a(x) F(u)dx.\end{equation}
 Denote by $q$ the conjugate of $p,$ i.e. $1/p+1/q=1.$ Recall the set $K_m$, given  in Definition \ref{def_cone}, and define $\psi: V \to (-\infty , \infty],$ by
\begin{equation}\label{gnov2}
\psi(u)= \left\{\begin{array}{ll}
\int_{\Omega} a(x) F^*(\frac{-\Delta u + u}{a(x)})\, dx, &   u\in K_m\cap W^{2,q}(\Omega) \\
+\infty,    &   u\notin K_m \cap W^{2,q}(\Omega), \\
\end {array}\right.
\end{equation}
with $Dom(\psi)=\{u \in V; \, \psi(u)< \infty\}.$  In Lemma \ref{lwr}, we shall  show that $\psi$ is convex and weakly  lower semi-continuous.\\

The following result is a direct consequence of Theorem \ref{con2}. However, for the convenience of the reader, we shall also  prove it in this paper.
\begin{corollary}\label{con3}  Assume that  $ u$ is a critical point of
\begin{equation}\label{eq2}I(w):= \psi(w)-\phi(w),
\end{equation}
where $\psi$ and $\phi$ are given in (\ref{gnov2}) and (\ref{gnov}) respectively.
If  there exists $v \in Dom(\psi)$ satisfying the linear equation,
\begin{equation} \label{need_100}
\left\{\begin{array}{ll}
-\Delta v + v= a(x)f(u), &  \mbox{ in } \Omega, \\
\frac{\partial v}{\partial \nu}= 0, &    \mbox{ on } \pOm,
\end {array}\right.
\end{equation}
then $u $ is a solution of the equation
\begin{equation*}
\left\{\begin{array}{ll}
-\Delta u + u= a(x)f(u), &   \mbox{ in } \Omega, \\
\frac{\partial u}{\partial \nu}= 0, &    \mbox{ on } \pOm.
\end {array}\right.
\end{equation*}
\end{corollary}
\noindent
\textbf{Proof.}
 Since $u$ is a critical point of $I,$ it follows from  Definition \ref{crit100} that
\begin{equation}\label{con33}\psi(w)-\psi(u) \geq \langle D \phi(u), w-u \rangle, \qquad \forall w\in V.\end{equation}
Since $I(u)$ is  finite  we have that $u \in Dom(\psi)$ and \[\psi(u)=\int_{\Omega} a(x) F^*\big (\frac{-\Delta u + u}{a(x)}\big )\, dx<\infty.\]
It then follows that $Au=-\Delta u + u \in L^q(\Omega)$ and $\psi(u)=\phi^*(A u)$ as shown in Lemma \ref{phi}.
By assumption, there exists $v \in Dom(\psi)$ satisfying $ A v= D \phi(u)$. Substituting $w=v$ in (\ref{con33}) yields that
\begin{equation}\label{con4}\phi^*(A v)-\phi^*(A u)=\psi(v)-\psi(u) \geq \langle D \phi(u), v-u \rangle=\langle A v, v-u \rangle. \end{equation}
On the other hand it follows from $ A v= D \phi(u)$ and Proposition \ref{var-pro} that $u \in \partial \phi^*(A v)$ from which we obtain
\begin{equation}\label{con5}\phi^*(A u)-\phi^*(A v) \geq \langle u, A u-A v \rangle. \end{equation}
Adding up (\ref{con4}) with (\ref{con5}) we obtain
\[\langle u, A u-A v \rangle+\langle A v, v-u \rangle\leq 0.\]
Since $A$ is symmetric we obtain that $\langle u-v, A u-Av \rangle \leq 0$ from which we obtain
\[\int_{\Omega}|\nabla u-\nabla v|^2 \, dx +\int_{\Omega}| u- v|^2 \, dx\leq 0,\]
thereby giving that $u=v.$ It then follows that $Au=Av=D \phi(u)$ as claimed.
\hfill$\square$\\

 Evidently, Corollary \ref{con3} maps out the plan for the prove of Theorem \ref{exi}. Indeed, by using Theorem \ref{MPT}, we show that the functional $I$ defined in (\ref{eq2}) has a nontrivial critical point and then we shall prove that the linear equation $(\ref{need_100})$ has a solution.\\

We shall need some preliminary results before proving the main Theorems in Introduction.
We first  list some properties of the function $F.$

\begin{lemma}\label{prop} The following assertions hold:
\begin{enumerate}
\item $F \in C^2(\R).$
\item $F: \R \to \R$ is strictly convex.
\item There exists a constant $C>0$ such that $0\leq F(t) \leq C(1+|t|^p)$ and
 $\mu F(t) \leq t F'(t)$ for all $t \in \R.$
 \item There exists  a constant $C>0$ such that  $F(t)\geq C t^\mu$ and $F'(t)t \geq C t^\mu$ for $|t|\geq 1.$
\item $F^* \in C^1(\R)$ and \[\frac{\mu}{\mu-1} F^*(s) \geq s {F^*}'(s), \qquad \forall s \in \R.\]
\item $F^*(s)\geq 0$ for all $s \in \R$ and $F^*(0)=0.$
\item There exists a constant $L>0$ such that $F^*(2s)\leq L F^*(s)$ for all $s \in \R.$
\end{enumerate}
\end{lemma}
\noindent
 \textbf{Proof.} Part 1) simply follows from $A_1.$  Part 2) is an immediate consequence of the fact $F^{''}(t)=f'(|t|) > 0$ for all $t\not=0$ and $F^{''}(0)=f'(0)=0.$ Part 3) follows from $A_2$ and $A_3.$  Part 4) follows from  part 3) and $A_1.$ We now proof part 5). The fact that  $F^* \in C^1(\R)$ follows from Theorem  26.6 in \cite{Roc}. Take $s \in \R$ and let $t_0$ be a point that maximizes $sup_{t \in \R} \{ts-F(t)\}$. Thus $s=F'(t_0).$ Since $F'$ is strictly monotone and hence invertible we have that $t_0 =(F')^{-1}(s)={F^*}'(s).$  It now follows from  $\mu F(t_0) \leq t_0 F'(t_0)$ and $F^*(s)=t_0s-F(t_0)$ that
\begin{eqnarray*}
 s {F^*}'(s)=st_0=t_0 F'(t_0)\geq \mu F(t_0)=\mu t_0s -\mu F^*(s)=\mu s {F^*}'(s)-\mu F^*(s)
 \end{eqnarray*}
 from which we obtain
 \[\mu F^*(s) \geq (\mu-1)s {F^*}'(s),\]
 as desired.\\
 Part 6): Since \[-F(0)=\inf_{s\in \R} F^*(s),\]
 we obtain that $F^*(s)\geq 0$ for all $s \in \R$. Also as  $0=f(0)= F'(0)=\partial F(0),$ it follows from Proposition \ref{var-pro} that \[F(0)+F^*(0)=0.\]
Thus, $F^*(0)=0$ as by the definition  $F(0)=0.$\\
Part 7):  By $A_3$, there exists a constant $l>0$ such that $2lF(t)\leq F(lt)$ for all $t \in \R.$   In the context of Orlicz spaces this property is known as $\nabla_2$ condition (see \cite{RR} for more details). It now  follows from $\nabla_2$ condition that for each $s\in \R,$
\begin{eqnarray*}
 F^*(s)=\sup_{t\in \R}\big \{st-F(t)\big \}\geq \sup_{t\in \R}\Big\{st-\frac{F(lt)}{2l} \Big\} &= & \frac{1}{2l}\sup_{t\in \R}\big\{2slt-F(lt)\big\}\\ &=& \frac{1}{2l}\sup_{t\in \R}\big\{2st-F(t)\big\}=\frac{1}{2l}F^*(2s).
 \end{eqnarray*}
 \hfill $\square$\\

Recall that $q$ is  the conjugate of $p,$ i.e. $\frac{1}{q}+\frac{1}{p}=1.$
\begin{lemma}\label{phi}
Assume that  $\phi: V \to \R$ is defined by $\phi(v)=\int_{\Omega} a(x) F(v)\, dx.$ Let $\phi^*: V^* \to (0, +\infty]$ be the Fenchel dual of $\phi.$  The following assertions hold.
\begin{enumerate}
\item
For each $h \in L^q(\Omega)$ we have
\[\phi^*(h)=\int a(x) F^*\big (\frac{h(x)}{a(x)}\big )\, dx.\]
\item  There exist positive constants $C_1$ and $C_2$ such that
\[\phi^*(h) \geq C_1\|h\|^q_{L^q(\Omega)}-C_2\]
for all $h \in L^q(\Omega).$
\item The function $\phi$ is differentiable and $\langle D \phi(u), u \rangle \geq \mu \phi(u)$ for all $u \in V.$
\item Let  $h \in  Dom (\phi^*)$. Then  $(1+t)h \in Dom (\phi^*)$ for all $0\leq t \leq 1.$ Moreover,  the directional derivative
\[D_{ h}\phi^*(h):=\lim_{t \to 0^+}\frac{\phi^*(h+t h)-\phi^*(h)}{t},\]
exists  and
\[0\leq D_{ h}\phi^*(h)= \int_\Omega {F^*}(\frac{h}{a})h \, dx\leq \frac{\mu}{\mu-1} \phi^*(h).
\]
\end{enumerate}
\end{lemma}
\noindent
\textbf{Proof.}
{\it 1.} Take $h\in L^q(\Omega).$ It follows from the density of $V$ in $L^p(\Omega)$ that
\begin{eqnarray*}
\phi^*(h)&=&\sup_{v \in V}\left\{\langle v, h \rangle -\phi(v) \right\}\\
&=&\sup_{v \in V}\left\{\int_\Omega v(x) h(x) \, dx -\int_\Omega a(|x|) F(v) \right\}\\
&=&\sup_{v \in L^p(\Omega)}\left\{\int_\Omega v(x) h(x) \, dx -\int_\Omega a(|x|) F(v) \right\}=\int_\Omega a(x) F^*\big (\frac{h(x)}{a(x)}\big )\, dx,
\end{eqnarray*}
where for the last equality we have used Proposition 2.1 in (\cite{Ek-Te}, page 271) and the fact that $a(x)>0$ for a.e. $x \in \Omega$.\\
{\it 2.} It follows from the boundedness of the function $a$ and part 3) of Lemma \ref{prop} that
\[\phi(v)=\int_{\Omega} a(x) F(v)\, dx \leq C \int (1+|v|^p)\, dx\]
for some constant $C>0$ and all $v\in L^p(\Omega).$ It then follows that
\begin{eqnarray*}
\phi^*(h)&=&\sup_{v \in V}\left\{\langle v, h \rangle -\phi(v) \right\}\\& \geq &
\sup_{v \in V}\left\{\int_\Omega v(x) h(x) \, dx  - C \int_\Omega (1+|v|^p)\, dx\right\}\\
&=&\sup_{v \in L^p(\Omega) }\left\{\int_\Omega v(x) h(x) \, dx  - C \int_\Omega (1+|v|^p)\, dx\right\}\\
&=& C_1\|h\|^q_{L^q(\Omega)}-C_2
\end{eqnarray*}
for some  constants $C_1$ and $C_2.$\\
{\it 3.} Differentiability of $\phi$ simply follows from $A_2$ and the fact that $a \in L^\infty (\Omega).$ An easy computation also shows that $D \phi(u)=a(x) F'(u).$ It now follows from part 3) of Lemma \ref{prop} that
\[\langle D \phi(u), u \rangle=\int_\Omega a(x)F'(u)u \, dx\geq \int_\Omega \mu a(x) F(u)\, dx=\mu \phi(u).\]
{\it 4.} It follows from part 7) of Lemma \ref{prop} that
there exists a constant $L>0$ such that $F^*$ satisfies the following condition,
\begin{equation}\label{rr}
F^*(2s) \leq L F^*(s), \qquad \forall s \in \R.
\end{equation}
Therefore, if $h \in Dom(\phi^*)$ then $\phi^*(2h)\leq L \phi^*(h) <\infty.$  On the other hand for   each $0<t<1$ by the convexity of $\phi^*$ we have that

\[ \phi^*((1+t)h)=\phi^*\big ((1-t)h+ t2h \big )\leq (1-t)\phi^*(h)+ t\phi^*(2h) < \infty,\]
from which we have that  $(1+t)h \in Dom(\phi^*).$\\

It  follows from  part 6) of Lemma \ref{prop}  that
 $ F^*(0)=\inf_{s \in \R}F^*(s).$ Since $F^*$ is differentiable we must have ${F^*}'(0)=0.$ For all  $s \in \R$ it follows from the monotonicity of ${F^*}'$ that

\[( {F^*}'(s)-{F^*}'(0))(s -0) \geq 0. \]
Thus,
\begin{eqnarray}\label{bbc0}s{F^*}'(s)\geq 0.\end{eqnarray}
It now follows from the latter inequality and part 5) of Lemma \ref{prop} that for all  $x \in \Omega$ with $a(x)\not =0$ and $h(x) \in \R,$  we have
\[0\leq 2h(x){F^*}'\big (\frac{2h(x)}{a(x)}\big )\leq \frac{\mu}{\mu-1} a(x)F^*\big (\frac{2h(x)}{a(x)}\big ).\]
Integrating both sides yields that
\[0\leq\int_\Omega 2h(x){F^*}'\big (\frac{2h(x)}{a(x)}\big )\, dx \leq  \frac{\mu}{\mu-1}\int_\Omega a(x)F^*\big (\frac{2h(x)}{a(x)}\big ) \, dx=\frac{\mu}{\mu-1}\phi^*(2h) <\infty.\]
It now follows from the monotonicity of ${F^*}'$ and the latter inequality that
\begin{eqnarray}\label{bbc1}
0\leq\int_\Omega h(x){F^*}'\big (\frac{(1+t)h(x)}{a(x)}\big )\, dx \leq \int_\Omega h(x){F^*}'\big (\frac{2h(x)}{a(x)}\big )\, dx \leq  \frac{\mu}{2(\mu-1)}\phi^*(2h)
\end{eqnarray}
By  the convexity of $F^*$ we have that the map
\[t \to \frac{ aF^*\big (\frac{h+ t  h}{a}\big) -a F^*\big (\frac{h}{a}\big )}{t}\]
is increasing  on $(0,1)$ and from  (\ref{bbc1}) we obtain  that
 \[0\leq  h(x){F^*}'\big (\frac{h(x)}{a(x)}\big )\leq \frac{ aF^*\big (\frac{h+ t  h}{a}\big) -a F^*\big (\frac{h}{a}\big )}{t}\leq h(x){F^*}'\big (\frac{(1+t)h(x)}{a(x)}\big )\leq h(x){F^*}'\big (\frac{2h(x)}{a(x)}\big ).\] It now follows from the dominated  convergence theorem and (\ref{bbc1}) that
 \begin{eqnarray*}\lim_{t \to 0^+}\frac{\phi^*(h+t h)-\phi^*(h)}{t} &=&\int_\Omega\lim_{t \to 0^+}\frac{ aF^*\big (\frac{h+ t  h}{a}\big) -a F^*\big (\frac{h}{a}\big )}{t}\, dx \\&=&\int_\Omega a{F^*}' (\frac{h}{a})  \frac{h}{a}\, dx\\&\leq& \int_\Omega \frac{\mu}{\mu-1}aF^*(\frac{h}{a})\, dx=\frac{\mu}{\mu-1} \phi^*(h),\end{eqnarray*}
from which the desired result follows.  \hfill $\square$\\

\begin{lemma}\label{lwr}
The functional  $\psi: V \to (-\infty , \infty]$ defined by
\begin{equation*}
\psi(u)= \left\{\begin{array}{ll}
\int_{\Omega} a(x) F^*(\frac{-\Delta u + u}{a(x)})\, dx, &   u\in K_m\cap W^{2,q}(\Omega) \\
+\infty,    &   u\notin K_m \cap W^{2,q}(\Omega), \\
\end {array}\right.
\end{equation*}
is convex and weakly lower semi-continuous.
\end{lemma}
\noindent
\textbf{Proof.} We first show that $\psi$ is  weakly lower semi-continuous. Let $\{u_n\}$ be a sequence in $V$ that converges weakly  to some $u \in V.$
If $\alpha:=\liminf_{n\to \infty}\psi(u_n)=\infty$ the there is nothing to prove. Let us assume that $\alpha<\infty$.  Thus, up to a subsequence, $u_n \to u$ a.e.,  $\psi(u_n)< \infty$ and $\lim_{n\to \infty}\psi(u_n)=\alpha.$ Since $u_n \to u$ a.e. we have that $u \in K_m.$ It follows from part 1) of Lemma \ref{phi} that $\psi(u_n)=\phi^*(-\Delta u_n + u_n)$. It also follows from part 2) of Lemma \ref{phi} that $\{u_n\}$ is also bounded in $W^{2,q}(\Omega).$ Thus, up to a subsequence,  we must have that  $u_n \to u$ weakly in $W^{2,q}(\Omega)$ and therefore  $u \in K_m \cap W^{2,q}(\Omega).$ Take $v \in L^p(\Omega).$ It follows that
\begin{eqnarray*}
\psi(u_n)=\phi^*(-\Delta u_n + u_n)\geq \int_\Omega v(x)(-\Delta u_n + u_n) \, dx-\phi(v),
\end{eqnarray*}
from which we obtain
\begin{eqnarray*}
\liminf_{n\to \infty}\psi(u_n)=\liminf_{n\to \infty}\phi^*(-\Delta u_n + u_n)\geq \int_\Omega v(x)(-\Delta u + u) \, dx-\phi(v),
\end{eqnarray*}
Taking sup over all $v \in L^p(\Omega)$ implies that
\begin{eqnarray*}
\liminf_{n\to \infty}\psi(u_n)=\liminf_{n\to \infty}\phi^*(-\Delta u_n + u_n)\geq \phi^*(-\Delta u + u)=\psi(u),
\end{eqnarray*}
from which the lower semi-continuity of $\psi$ follows.\\

We now show that $\psi$ is convex. Let $u_1 ,u_2 \in V$ and $t \in (0,1).$ We need to verify that
\begin{equation}\label{waycon}\psi\big (t u_1+(1-t)u_2\big ) \leq t \psi(u_1)+(1-t)\psi(u_2).\end{equation}
Note first that  $F^*$ is  non-negative by Lemma \ref{prop}, and  $a\geq 0$ by assumption. Thus,  we have that $\psi \geq 0.$  If one of $\psi(u_1)$ or $\psi (u_2)$ is $+\infty$ then we are done. So assume that $\psi(u_1), \psi (u_2) \in \R.$  It then follows that $u_1, u_2 \in K_m\cap W^{2,q}(\Omega)$ and
\begin{equation}
\int_{\Omega} a(x) F^*\Big(\frac{-\Delta u_i + u_i}{a(x)}\Big)\, dx < \infty, \qquad i=1,2.
\end{equation}
Since $K_m\cap W^{2,q}(\Omega)$ is a convex set we have that $t u_1+(1-t)u_2 \in K_m\cap W^{2,q}(\Omega)$. On the other hand,
for almost every $x \in \Omega,$ it follows from the convexity of $F^*$ and linearity of the map  $u \to -\Delta u+ u$ that
\begin{eqnarray*}
F^*\Big(\frac{-\Delta (t u_1+(1-t)u_2) + t u_1+(1-t)u_2}{a(x)}\Big)&=&F^*\Big(\frac{t(-\Delta  u_1+u_1)+(1-t)(-\Delta u_2+u_2)}{a(x)}\Big)\\ &\leq &
tF^*\Big(\frac{-\Delta  u_1+u_1}{a(x)}\Big)+(1-t)F^*\Big(\frac{-\Delta u_2+u_2}{a(x)}\Big).
\end{eqnarray*}
Therefore,  by multiplying the latter expression by $a(x)$ and integrating over $\Omega$  the inequity (\ref{waycon}) follows.
 \hfill $\square$\\

\begin{lemma}\label{nn}
Let $\{u_n\} \subset K_m$ be a sequence in $H^1(\Omega)$ that converges weakly to some $u \in K_m$ and also $u_n \to u$ a.e.. Then
\begin{enumerate}
\item  For each $2\leq r< 2_m^*$ we have that $u_n \to u$ strongly in $L^r(\Omega).$
\item  $\phi(u_n) \to \phi(u).$
\item $\langle D \phi(u_n), u_n-u \rangle \to 0.$
\item For each $\delta>0,$ we have $\phi(\delta u_n-\delta u) \to 0.$
\end{enumerate}
\end{lemma}
\noindent
\textbf{Proof.} 1)  Suppose $ 2 \le r < 2_m^*$ and set $ T:=\min\{ 2_m^*,r+1 \}$ (recall if $m=2$ then $ 2_m^*=\infty$).  By interpolating $L^2$ and $L^T$ we have

\[\|u_n-u\|_{L^r}\leq \|u_n-u\|^\theta_{L^2}\|u_n-u\|^{1-\theta}_{L^{T}},\]
for some $0<\theta \leq 1.$ Since $u_n, u \in K$, it follows from Lemma \ref{imbeddings} that $\|u_n-u\|_{L^{T}}$ is bounded. It now follows from the weak convergence in $H^1$ that $u_n \to u$ strongly in $L^2,$ from which together with latter inequality we get that $u_n \to u$ strongly in $L^r(\Omega).$

2)  It follows from $A_2$ that $F(u_n) \leq C(1+|u_n|^p).$ Thus, the result follows from part 1) and the dominated convergence theorem.\\
3) Note that
\[|\langle D \phi(u_n), u_n-u \rangle|\leq \int a(|x|)|f(u_n)( u_n-u)|\,  dx \leq C \int |u_n-u|(1+|u_n|^{p-1})\, dx,\]
 and by Holder inequality and the result of part 1) we obtain
\[\int |u_n-u|(1+|u_n|^{p-1})\, dx\leq \int |u_n-u|(1+|u_n|)^{p-1}\, dx\leq \|u_n-u\|_{L^p}\|1+|u_n|\|_{L^p}^{\frac{p}{q}} \to 0.\]
 4)
 It follows from $A_2$ that $F(\delta u_n- \delta u ) \leq C(1+|\delta u_n- \delta u|^p)$, which together with the dominated convergence theorem and part 1) the desired result follows.
\hfill $\square$\\

\begin{proposition}\label{critical}
Consider  the functional $I: V\rightarrow \mathbb{R}\cup \{+\infty\}$ defined by  $$I(u):= \psi(u)- \varphi(u),$$
where  $\phi(u)=\int _{\Omega} a(x) F(u) \, dx$ and
 \begin{equation*}
\psi(u)= \left\{\begin{array}{ll}
\int_{\Omega} a(x) F^*(\frac{-\Delta u + u}{a(x)})\, dx, &   u\in K_m\cap W^{2,q}(\Omega) \\
+\infty,    &   u\notin K_m\cap W^{2,q}(\Omega). \\
\end {array}\right.
\end{equation*}

Then $I$ has a nontrivial  critical point.
\end{proposition}
\noindent
\textbf{Proof.} We make  use Theorem \ref{MPT} to prove this lemma.
 First note  that, by $A_2$, the functional $\phi$ is $C^1$
and
\[D\varphi(u)=  a(x)f(u).\]
 Note also that
 $\psi$ is proper and convex  as $K_m \cap W^{2,q}(\Omega)$ is convex  in $V.$ It also follows from Lemma \ref{lwr} that $\psi$ is weakly lower semi-continuous. We shall now proceed  in several steps. \\

\noindent
{\it Step 1.} In this step we shall verify the mountain pass geometry  for $I$.\\
By Lemma \ref{prop} we have that  $\mu F(t) \leq t F'(t)$ and $\frac{\mu}{\mu-1} F^*(t) \geq t {F^*}'(t)$. Thus there exist constants $C_1$ and $C_2$ such that for $t \geq 1,$
\begin{equation}\label{ini}
F(t) \geq C_1 |t|^\mu, \quad F^*(t)\leq C_2 |t|^{\frac{\mu}{\mu-1}}.
\end{equation}

It is clear that $I(0)=0$ as $F^*(0)$ by part 6) of Lemma \ref{prop}. Since  $e=\sup_{x \in \Omega} a(x)+1 \in K_m$,  it follows that for $t \geq 1$
\begin{eqnarray*}\label{MPG1}
 I(te)&=&
\int a(x)F^*\big (\frac{te}{a(x)}\big )dx- \int_{B_1} a(|x|)F(te) dx\\
 &\leq & C_2 \int a(x)^{1-\frac {\mu}{\mu-1}}|t e |^{\frac {\mu}{\mu-1}} dx- C_1\int a(|x|)|te|^{\mu} dx\\
 \end{eqnarray*}
Now, since $\mu>2$ one has that $\frac {\mu}{\mu-1}<2.$ Thus for $t$ sufficiently large $I(te)$ is negative. We now prove condition 3) of $(MPG)$. Take
 $u\in Dom(\psi)$ with $\|u\|_{V}= \rho>0$. We have
\begin{equation}\label{MPG6}
  I(u)=\psi(u)-\phi(u)= \varphi^*(A u)- \varphi(u)\geq <A u,u>- 2\varphi(u)= \|u\|^2_{H^1}- 2\varphi(u)
\end{equation}
Since the function $a(x)$ is bounded it follows from $A_2$ that
\begin{equation}\label{grow}
\forall \delta>0, \, \exists C_\delta>0, \qquad |a(x)F(t)| \leq \delta |t|^2+ C_\delta |t|^p,
\end{equation}
from which we obtain
\begin{equation}\label{grow1}
\forall \delta>0, \, \exists C_\delta>0, \qquad \phi(u)=\int a(x)F(u) \leq \delta \|u\|^2_{L^2(\Omega)}+ C_\delta\|u\|^p_{L^p(\Omega)}.
\end{equation}
Note that from Lemma \ref{imbeddings}, for $u\in K_m$ one has
$\|u\|_{L^p}\leq C_p \|u\|_{H^1}$ and therefore,
\begin{equation}\label{MPG7}
 \|u\|_V= \|u\|_{H^1}+ \|u\|_{L^p}\leq (1+ C_p)\|u\|_{H^1}
\end{equation}
for some constant $C>0.$
It now follows from (\ref{grow1}) and (\ref{MPG7})  that
\begin{eqnarray*}I(u)&=& \|u\|^2_{H^1}- 2\varphi(u)\geq \|u\|^2_{H^1} -2 \delta \|u\|^2_{L^2(\Omega)}-2 C_\delta \|u\|^p_{L^p(\Omega)}\\
&\geq &(1-2 \delta) \|u\|^2_{H^1} -2 C_\delta C^p_p \|u\|^p_{H^1(\Omega)}\\
& \geq &\frac{(1-2 \delta) }{(1+C_p)^2}\|u\|^2_{V}-2 C_\delta C^p_p \|u\|^p_{V}\\
&=&\frac{(1-2 \delta) }{(1+C_p)^2}\rho^2-2 C_\delta C^p_p\rho^p
\end{eqnarray*}
Therefore, for $\delta <1/2$  we have that

\begin{equation*}
 I[u]\geq \frac{(1-2 \delta) }{(1+C_p)^2}\rho^2-2 C_\delta C^p_p\rho^p>0
\end{equation*}
provided  $\rho>0$ is  small enough. If $u\notin Dom(\psi)$, then clearly $I(u)>0$.  Therefore (MPG) holds for the functional $I$.\\

\noindent
{\it Step 2.} We verify Palais-Smale  compactness condition.
Suppose that $\{u_n\}$ is a sequence in $K_m$ such that $I(u_n)\rightarrow c\in \mathbb{R}$ as $\epsilon_n \rightarrow 0$ and
\begin{equation}\label{d}
\langle D\varphi(u_n), u_n-v \rangle + \psi(Av)- \psi(Au_n)\geq -\epsilon_n \|v- u_n\|_V, \qquad \forall v\in V.
\end{equation}
We must show that $\{u_n\}$ has a convergent subsequence in $V$. First, note that  $u_n\in Dom(\psi)$ and therefore,   \[I(u_n)=\varphi^*(Au_n)- \varphi(u_n)\rightarrow c, \quad \text{ as } n \to \infty.\]
Thus, for large values of $n$ we have
\begin{equation}\label{d0}
\varphi^*(Au_n)- \varphi(u_n)\leq 1+c.
\end{equation}

Since $Au_n \in Dom(\phi^*),$ it follows from part 4) of Lemma \ref{phi} that $(1+r)Au_n \in Dom(\phi^*)$ for  $0\leq  r\leq 1$.
By setting  $v= (1+r)u_n \in K_m \cap W^{2,q}(\Omega)$ for  $0< r \leq 1$ in  (\ref{d}) we have that
\begin{equation}\label{d1}
-\langle D\varphi(u_n),r u_n \rangle + \phi^*(Au_n+ r Au_n)- \phi^*(Au_n)\geq -r\epsilon_n \| u_n\|_V.
\end{equation}
Dividing both sides by $r$ and letting $r \to 0^+$ yield that,
\begin{equation}\label{d2}
-\langle D\varphi(u_n), u_n \rangle + D_{Au_n} \phi^*(Au_n) \geq -\epsilon_n \| u_n\|_V,
\end{equation}
where $D_{Au_n} \phi^*(Au_n) $ is the directional derivative of $\phi^*$ at $Au_n$ in the direction $Au_n$ that exists due to Lemma \ref{phi} part 4), and furthermore $D_{Au_n} \phi^*(Au_n)\leq \frac{\mu}{\mu-1}\phi^*(Au_n).$
Multiply (\ref{d2}) by $-1/2$ and sum it up  with  (\ref{d0}) to get
\begin{equation*}
\varphi^*(Au_n)-\frac{1}{2}D_{Au_n} \phi^*(Au_n)  +\frac{1}{2}\langle D\varphi(u_n), u_n \rangle -\varphi(u_n)\leq  1+c+ \| u_n\|_{V},
\end{equation*}
and therefore by 3) and 4) in Lemma \ref{phi} we obtain that
\begin{equation}\label{d4}
\big (1-\frac{\mu}{2(\mu-1)}\big )\varphi^*(Au_n)+ \big (\frac{\mu}{2}-1)\varphi(u_n)\leq  1+c+ \| u_n\|_{V}.
\end{equation}
Since $\mu>2$ we have that
\[\frac{\mu}{2}-1>0 \qquad \&\qquad   1-\frac{\mu}{2(\mu-1)}>0.\]
Taking now into account that $\varphi^*(Au_n)\geq 0, \varphi(u_n)\geq 0$, it follows from (\ref{d4}) that
\begin{equation} 
 \varphi^*(Au_n)+ \varphi(u_n)\leq C_2( 1+ \| u_n\|_{V}),
\end{equation}
for an appropriate constant $C_2>0$. On the other hand
$$ \varphi^*(Au_n)+ \varphi(u_n)\geq \langle Au_n, u_n\rangle = \| u_n\|_{H^1}^2,$$
which according to (\ref{d4}) results in
\begin{equation*}
\| u_n\|_{H^1}^2\leq C_2( 1+ \| u_n\|_{V}).
\end{equation*}
It also follows (\ref{MPG7}) that $\| u_n\|_{V} \leq (1+C_p)\| u_n\|_{H^1}$ and therefore
\[\| u_n\|_{H^1}^2\leq C_0( 1+ \| u_n\|_{H^1}),\]
for some constant $C_0.$
Therefore  $\{u_n\}$ is bounded in $H^1(\Omega).$ By  passing to a subsequence if necessary, there exists $\bar u \in H^1(\Omega)$ such that  $u_n\rightharpoonup \bar{u}$ weakly in $H^1(\Omega)$,  $u_n\rightarrow  \bar{u}$ strongly in $L^2(\Omega)$ and $u_n\rightarrow  \bar{u}$ a.e..
Note first that  $\bar{u}\in K_m$.
It also follows from (\ref{d4}) that $\{\phi^*(Au_n)\}$ is bounded and therefore,  by Lemma \ref{lwr}, we obtain
\[ \varphi^*(A\bar{u})\leq \liminf_{n\rightarrow \infty} \varphi^*(Au_n)< \infty,\]
from which one has  $\bar u \in Dom(\psi).$
By setting  $v= \bar{u}$ in (\ref{d}) we obtain
\begin{equation}\label{der1}
-\langle D \phi(u_n), \bar{u}- u_n\rangle+ \varphi^*(A\bar{u})- \varphi^*(Au_n)\geq -\epsilon_n \|\bar{u}- u_n\|_V, \end{equation}
By Lemma \ref{nn} we have that $\langle D \phi(u_n), \bar{u}- u_n\rangle \to 0,$ and by (\ref{MPG7}) we have that $\|u_n-\bar u\|_V$ is bounded.
Therefore passing into limits in (\ref{der1}) results in
\begin{equation}\label{der2}
\limsup_{n\rightarrow \infty} \varphi^*(Au_n)\leq \varphi^*(A\bar{u}).
 \end{equation}
The latter inequality together with the fact that $\varphi^*(A\bar{u})\leq \liminf_{n\rightarrow \infty} \varphi^*(Au_n)$ yield that
\[\varphi^*(A\bar{u})= \lim_{n\rightarrow \infty} \varphi^*(Au_n).\]
Now observe that
\begin{align}\label{conv0}
\|u_n\|_{H^1}^2- \|\bar{u}\|_{H^1}^2= \langle Au_n,u_n\rangle - \langle A{\bar{u}},\bar{u}\rangle
 =\langle Au_n,u_n- \bar{u}\rangle+ \langle Au_n- A\bar{u},\bar{u}\rangle.
\end{align}
But  weakly convergence of $u_n$ to $\bar{u}$ in $H^1(\Omega)$ means that
$Au_n \rightharpoonup A\bar{u}$ weakly in $H^{-1}(\Omega)$, thus
\begin{equation}\label{conv1}
\langle Au_n- A\bar{u},\bar{u}\rangle\rightarrow 0, \quad \text{ as } n \to \infty.
 \end{equation}
Let $0<\delta<1.$ It follows from  Lemma \ref{phi} and part 6) of Lemma \ref{prop} that $\phi^*(0)=0.$ Thus, by the convexity of $\phi^*$ we have that
\[\phi^*(\delta A u_n)=\phi^*\big (\delta A u_n +(1-\delta)0\big)\leq \delta \phi^*(A u_n)+(1-\delta)\phi^*(0)=\delta \phi^*(A u_n).\]
We then have that
\begin{eqnarray*}
|\langle Au_n,u_n- \bar{u}\rangle |&\leq & \phi(\frac{u_n-\bar u}{\delta})+\phi^*(\delta A u_n)\\
&\leq & \phi(\frac{u_n-\bar u}{\delta})+\delta \phi^*(A u_n).\\
\end{eqnarray*}
By taking $\limsup$ as $n \to \infty$ we have that
\[\limsup_{n \to \infty}|\langle Au_n,u_n- \bar{u}\rangle |\leq  \limsup_{n\to \infty }\phi(\frac{u_n-\bar u}{\delta})+\delta\limsup_{n\to \infty} \phi^*(A u_n)\]
Now by virtue of Lemma \ref{nn} we have that $\limsup_{n\to \infty }\phi(\frac{u_n-\bar u}{\delta})=0$ from which we obtain
\[\limsup_{n \to \infty}|\langle Au_n,u_n- \bar{u}\rangle \leq \delta \varphi^*(A\bar{u}).\]
By now letting $\delta \to 0$ we obtain  that

\begin{equation}\label{conv2}\limsup_{n \to \infty}|\langle Au_n,u_n- \bar{u}\rangle|=0.\end{equation}
Therefore, from (\ref{conv0}), (\ref{conv1}) and (\ref{conv2}) one has
$$u_n\rightarrow \bar{u}\quad \text{strongly in }\quad H^1.$$
It now follows from Lemma \ref{nn} part 1) that $u_n \to \bar u$ strongly in $L^p(\Omega).$ Therefore,
$$u_n\rightarrow \bar{u}\quad \text{strongly in }\quad V,$$
as desired.
\hfill $\square$\\

\noindent
\textbf{Proof of Theorem \ref{exi}.} It follows from Proposition \ref{critical} that
 the functional $I$ has a nontrivial critical point $u \in K_m$.   We will now apply
Corollary \ref{con3} to see that $u $ is nonnegative nontrivial   monotonic solution
of (\ref{eq}).  To do this we need to show there is some $v  \in Dom(\psi)$ which
satisfies $A v=D \phi(u)$; or to be more explicit,  which satisfies
(\ref{need_100}).   We now prove this.  Fix $ u \in K_m$ and suppose $a$ satisfies
the assumed hypothesis.  For $ \E>0$ small let $ u^\E,a^\E$ denote the smoothed
versions of $ u$ and $ a$ respectively as promised by Lemma \ref{smoothing}.
Replacing $ af(u)$ with $ a^\E f(u^\E)$ on the right hand side of (\ref{need_100})
we can apply Proposition \ref{mono_linear} to see there is some $ v^\E \in K_m \cap
C^{2,\alpha}(\overline{\Omega})$ which satisfies $ -\Delta v^\E + v^\E = a^\E
f(u^\E)$ in $ \Omega$ with $ \partial_\nu v^\E=0$ on $ \pOm$, which has the weak
formulation
\begin{equation} \label{weak_form_100}
\int_\Omega \nabla v^\E \cdot \nabla \phi +  v^\E \phi dx = \int_\Omega a^\E f(u^\E)
\phi dx, \quad \forall \phi \in H_G^1(\Omega).
\end{equation}  Note note that we have $ |f(u^\E)| \le C ( 1+ | u^\E|^{(p-1)}) \le C
(1+ |u|^{(p-1)})$ in $ \Omega$ and so we have $ |f(u^\E)|^{p'} \le C_p ( 1 + |
u^\E|^{(p-1)p'}) = C_p ( (1 + | u^\E|^p)$ in $ \Omega$ and hence  $
\|f(u^\E)\|_{L^{p'}}^{p'} \le C (1 + \|u^\E\|_{L^p}^p ) \le C (1+ \|u\|_{L^p}^p)$.

Taking $ \phi=v^\E$ in (\ref{weak_form_100}) and applying  H\"older's inequality on
the right hand side one obtains, after using the above bound,
 \[ \|v^\E \|_{H^1}^2 \le \|a^\E f(u^\E) \|_{L^{p'}} \|v^\E\|_{L^{p}} \le C \|v^\E
\|_{H^1}\]
 where $C$ independent of $ \E$ and where we have used the improved Sobolev
imbeddings for $v^\E \in K_m$.   This gives an $H^1(\Omega)$ bound on $ v^\E$ and
hence by passing to a suitable subsequence there is some $ v \in H^1_G(\Omega)$
such that $ v^\E \rightharpoonup v$ in $H_G^1(\Omega)$ as $ \E \rightarrow 0$.  We
can then pass to the limit in (\ref{weak_form_100}) for all $ \phi \in
H_G^1(\Omega) \cap L^\infty(\Omega)$;  to pass to the limit on the right hand side
we can use the dominated convergence theorem.  Noting that $K_m$ is weakly closed
in $H_G^1(\Omega)$ we have $ v \in K_m$.    To complete showing that $v \in
Dom(\psi)$ we need to show that $v \in W^{2,q}(\Omega)$ where $ q=p'$.    Note
above that we have shown $ a^\E f(u^\E)$ is bounded in $L^{p'}(\Omega)$
independently of $ \E$.  So apply $L^p$ elliptic regularity theory shows that $
v^\E$ is bounded in $ W^{2, {p'}}(\Omega)$.    So we now have a nonnegative nonzero sufficiently
regularity solution $u$ of (\ref{eq}) and we can then apply the strong maximum
principe to see that $u$ is positive.

 \hfill $\square$\\

\noindent
\textbf{Proof of Corollary \ref{exx}.}   Suppose $ \Omega$ is a domain of $m$ revolution which satisfies (\ref{assume_cube}) and    that $ a$ satisfies $A_4$ and further we assume that
$ a(t_1,t_2, ..., t_m)= a(t_1, t_2, ... , t_i)$  for some $ 1 \le i <m$.
Suppose  $A_1-A_4$ hold with $ 2< p<2_i^*$ in $A_2$.
Applying Theorem \ref{exi} one sees there is a positive  $ v=v(t_1, ... ,t_i) \in K_i$ which satisfies the lower dimensional problem
 \begin{equation} \label{nonnon_low}
 -\Delta_t v - \sum_{k=1}^i \frac{n_k-1}{t_k} v_{t_k} + v = a(t_1,..,t_i) f(v) \mbox{ in } Q_i, \qquad \mbox{$ \partial_\nu v=0$ on $ \partial Q_i$.}
 \end{equation}
Set  $ u(t_1, ... , t_m):= v(t_1, ... , t_i)$ and note that $u \in K_m$ is a nonzero solution
\begin{equation} \label{nonnon}
 -\Delta_t u - \sum_{k=1}^m \frac{n_k-1}{t_k} u_{t_k} + u = a(t_1,...,t_m) f(u) \mbox{ in } Q_m, \qquad \mbox{$ \partial_\nu u=0$ on $ \partial Q_m$.}
 \end{equation}
\hfill $\Box$

\section{Appendix}

\begin{lemma} \label{smoothing}  (Smoothing of $ u \in K_m$) Suppose $m \ge 2$,  $ Q_m=\tilde{\Omega}$  and  $ u \in K_m$.  Then there is some smooth $ u^\E \in K_m$ such that $ u^\E \le u$ a.e. in $\Omega$ and such that $u^\E \rightarrow u$ a.e. in $ \Omega$ as $ \E \searrow 0$.

\end{lemma}

\noindent
\textbf{Proof.} Consider $ 0 \le u \in K_m$ and so $0 \le u \in Y_m$ with $ u_{t_k} \ge 0$ a.e. in $ Q_m$ for $ 1 \le k \le m$.     Let $ 0 \le \eta$ denote a smooth compactly supported function in $(-1,0)^m$ with $ \int_{\IR^m} \eta =1$.    \\

Note that $u$ is defined in $ Q_m$ and we then extend $ u$ to be zero outside $Q_m$.  Let $ \E>0$ be small and  define $ Q^\E:= (-\E,0)^m$ and we let $ t=(t_1,t_2,..,t_m)$ and $ \overline{t}:=( \oo{t}_1,..., \oo{t}_m)$.  We then define
\[ u^\E(t):=\int_{Q^\E} \frac{1}{\E^m} \eta \left(  \frac{ \oo{t}}{\E} \right) u(t+ \oo{t}) d \oo{t},\] for $ t \in Q_m$.   Note that we can re-write $u^\E$ as
\[ u^\E(t)= \int_{\IR^m} \frac{1}{\E^m} \eta \left(  \frac{ \hat{t}-t}{\E} \right) u( \hat{t}) d \hat{t},\] and this shows that $ u^\E$ is smooth.   Let $ t \in Q_m$ and returning to the first expression we see
\[ u^\E(t) = \int_{Q^\E} \frac{1}{\E^m} \eta \left(  \frac{ \oo{t}}{\E} \right) u(t+ \oo{t}) d \oo{t} \le \int_{Q^\E} \frac{1}{\E^m} \eta \left(  \frac{ \oo{t}}{\E} \right) u(t) d \oo{t} =u(t),\] after recalling the support of $ \eta$ and since $u$ has some monotonicity.  Now suppose that $ t=(t_1,...,t_m),s=(s_1,...,s_m) \in Q_m$ with $ t_k \le s_k$ for all $ 1 \le k \le m$.   Then note
\[   u^\E(t) = \int_{Q^\E} \frac{1}{\E^m} \eta \left(  \frac{ \oo{t}}{\E} \right) u(t+ \oo{t}) d \oo{t} \le \int_{Q^\E} \frac{1}{\E^m} \eta \left(  \frac{ \oo{t}}{\E} \right) u(s+ \oo{t}) d \oo{t} = u^\E(s),\] and hence we see $u^\E$ has the desired monotonicity.   One can then use a the standard approach to show that $u^\E \rightarrow u$ a.e. in $ Q_m$ as $ \E \searrow 0$.

\hfill $\square$\\

\noindent
\textbf{Proof of Proposition \ref{mono_linear}.} Assume $ \Omega$ is domain of $m$ revolution in $ \IR^N$ which satisfies
(\ref{assume_cube}) and   fix $ 0 \le h \in K_m \cap C^{1,1}(\oo{\Omega})$.
    As before we identify $H_G^1(\Omega)$ and $Y_m$ and hence,  by standard
arguments there is a unique $ v \in Y_m$ which  satisfies
\begin{equation} \label{weak_u}
 \int_{Q_m}  \left( \nabla_t v \cdot \nabla_t \eta + v \eta \right) d \mu_m(t)
=\int_{Q_m} h(t) \eta \; d \mu_m(t), \quad \forall \eta \in Y_m.
\end{equation}

 Our goal now is to show that $ v \in K_m \cap
C^{2,\alpha}(\oo{\Omega})  $ and we first show that $v \in K_m$.    To do that we
begin by solving a smoothed version of (\ref{weak_u}). For $\E>0$ small define
$ d \mu^\E_m(t)= \prod_{i=1}^m (t_i+\E)^{n_i-1} d t_i$ and consider the energy
\[ E_\E(v):=\frac{1}{2} \int_{Q_m} \left( | \nabla_t v|^2 + v^2 - h(t) v \right) d
\mu_m^\E(t).\]  It is easily seen that $ E_\E$ attains its minimum on $H^1(Q_m)
\subset Y_m$ at $ v^\E$ which satisfies
 \begin{equation}  \label{weak_reg}
 \left\{
\begin{array}{rrl}
   -\Delta_t v^\E -\sum_{k=1}^m \frac{n_k-1}{t_k+\E} v_{t_k}^\E + v^\E    &=& h
\qquad \hfill \mbox{ in } Q_m, \\
\partial_\nu v^\E&=& 0 \qquad \hfill  \mbox{ on } \partial Q_m.
\end{array}
\right.
\end{equation}
 Note now that   the solution of this problem is as smooth as the right hand side and the nonsmooth domain allow. We now proceed in several steps.
   \\

\noindent
 {\it Step 1.} For $0<\alpha<1$ and $ 0<\E <1$ one has  $ v^\E \in
C^{2,\alpha}(\overline{Q_m})$.   We prove this result in the case of a double
domain of revolution;  but it extends to the general case.   We will use the method
 of even reflections to prove the global regularity result. Towards this define  $
c^\E_k(t)=\frac{n_k-1}{t_k+\E}$ and let $ c^\E(t)$ be the $2$ dimensional vector
with components  $c^\E_1(t)$ and $c^\E_2(t)$.  Define the even extension of $v^\E$ by
\begin{equation*}
 \overline{v^\E}(t_1,t_2):= \left\{
\begin{array}{rrl}
v^\E(t_1,t_2) && \qquad \mbox{ in } Q_2 \\
v^\E(-t_1,t_2)   & &  \qquad    \mbox{ in } (-1,0) \times (0,1)  \\
v^\E(-t_1,-t_2)  &&   \qquad \mbox{ in } (-1,0) \times (-1,0) \\
v^\E(t_1,-t_t) &&  \qquad \mbox{ in } (0,1) \times (-1,0)
\end{array}
\right.
\end{equation*} and where $\oo{v^\E}$ is extended to the axis by continuity.    We
now define $\oo{c^\E_k}(t)$ to be the odd reflection of $c^\E_k$ across $ t_k=0$,
ie. given by
\begin{equation*}
 \overline{c^\E_k}(t):= \left\{
\begin{array}{rrl}
c^\E_k(t) && \qquad \mbox{ in } t_k \ge 0\\
-c^\E_k(-t)  & &  \qquad    \mbox{ in } t_k<0, \\
\end{array}
\right.
\end{equation*} and note that $\oo{c^\E_k}$ has a jump discontinuity at $t_k=0$.
Set $ \oo{c^\E}$ to be the $2$ dimensional vector with components $ \oo{c^\E_k}(t)$.

We then let $\oo{h}$ denote the even extension of $h$  as we did with $ v^\E$.
Consider the nested cubes
\[  D_4 \subset D_3 \subset D_2
\subset D_1 \subset  \left( -1, 1 \right)\times \left( -1, 1 \right),\] where each is
compactly contained in the other.
We then have $\oo{v^\E}$ is a weak solution of
\begin{equation} \label{loppe}
-\Delta_t \oo{v^\E} + \oo{v^\E} = \oo{h} + \oo{c^\E}(t) \cdot \nabla_t \oo{v^\E}
\quad \mbox{ in } D_1.
  \end{equation}  Note that since $ v^\E \in H^1(Q_2)$ we have $ \oo{v^\E} \in
H^1(D_1)$ which implies the right hand side of (\ref{loppe}) belongs to
$L^2(D_1)$.

     We can then apply elliptic regularity to see that $ \oo{v} \in H^2(D_2)$ and
since we are in dimension 2 we can apply the Sobolev imbedding theorem to see
that  $\oo{v^\E}_{t_k} \in L^q(D_2)$ for all $ 1\leq q<\infty$.  In the case of $m>2$
one can apply a bootstrap argument to obtain $\oo{v^\E}_{t_k} \in L^q(D_2)$ for
all $ 1\leq  q< \infty$.
  We can then apply $L^q$ elliptic regularity theory and the Sobolev imbedding
theorem to see that $ \oo{v^\E} \in C^{1,\alpha}(D_3)$ for all $ 0<\alpha<1$.
We now show the right hand side of (\ref{loppe}) is H\"older continuous on
$D_3$.\\

\noindent {\it Claim 1.}  $  (\oo{c^\E} \cdot \nabla_t \oo{v^\E})  \in
C^{0,\alpha}(D_3)$.  This is not immediately obvious since $ \oo{c^\E}_k$ has jump
discontinuities.   We now show that $ \oo{c^\E}_1 \oo{v^\E}_{t_1} \in
C^{0,\alpha}(D_3)$. \\
Consider $ (t_1,t_2), ( \tau_1, \tau_2) \in D_3$ and consider the three cases: \\
 \noindent
 (i) $ t_1 <0$ and $ \tau_1 \ge 0$, \quad
(ii) $t_1,\tau_1 \ge 0$, \quad
(iii) $ t_1,\tau_1<0$. \\

\noindent

\noindent
Case (i).  Now note that
\[
\big| ( \oo{c^\E}_1 \oo{v^\E}_{t_1})(t_1,t_2) - ( \oo{c^\E}_1
\oo{v^\E}_{t_1})(\tau_1,\tau_2)\big| \le C \left(  \big|  \oo{v}_{t_1}(t_1,t_2)\big|
+ \big| \oo{v^\E}_{t_1}(\tau_1,\tau_2)\big| \right),\] and using the the fact that $
v^\E_{t_1}=0$ on $ t_1=0$ and $ \oo{v^\E}_{t_1} \in C^{0,\alpha}(D_3)$ we have
\[ \big|  \oo{v^\E}_{t_1}(t_1,t_2)\big| = \big|  \oo{v^\E}_{t_1}(t_1,t_2) -
\oo{v^\E}_{t_1}(0,t_2) \big| \le C | t_1|^\alpha,\] and similarly we have $|
\oo{v^\E}_{t_1}(\tau_1,\tau_2) | \le C |\tau_1|^\alpha$.  From this we see
\[
\frac{\big| ( \oo{c^\E}_1 \oo{v^\E}_{t_1})(t_1,t_2) - ( \oo{c}_1
\oo{v^\E}_{t_1})(\tau_1,\tau_2)\big|}{ \big| (t_1,t_2)-(\tau_1,\tau_2) \big|^\alpha}
\le \frac{ C |t_1|^\alpha + C |\tau_1|^\alpha }{ \big| ( t_1 -\tau_1, t_2-\tau_2)
\big|^\alpha},\] and note the right hand side is bounded after recalling that $ t_1
<0$ and $ \tau_1 \ge 0$.   \\

For case (ii) and (iii) we easily see the needed H\"older quotient is bounded after
we consider that $\oo{c^\E_1}$ is smooth on the restricted domain and since $
\oo{v^\E}_{t_1}$ is H\"older continuous.   This completes the proof of the claim.
 \\

Using the above claim we now see that the right hand side of (\ref{loppe}) is in
$C^{0,\alpha}( D_3)$ and hence we can now apply elliptic regularity theory to see
that $ \oo{v^\E} \in C^{2,\alpha}( \overline{D_4})$.      We now argue that one in
fact has $ v^\E \in C^{2,\alpha}(\overline{Q_2})$.        We could extend $ v^\E$
evenly across the outer boundaries and the extension would satisfy a similar
equation as above.  Note one difference now is that $\oo{h}$ is sufficiently smooth
and  symmetric across $ t_k=0$ and hence its  even extension is sufficiently
regular.  When we extend $h$ across the outer boundary it will, in general, only be
Lipschitz continuous.   But this is enough to carry out the above procedure.
Carrying this out gives $ v^\E \in C^{2,\alpha}(\overline{Q_2})$, which complete the
proof of Step  1.   \\

\noindent
 {\it Step 2.} We now show that $ v^\E \in K_m$.    We do this step for all $m\geq 2$.   First
note that we have $ v^\E \ge 0$ and so we need only show that $v^\E$ has the
desired monotonicity; $ v^\E_{t_k } \ge 0$ in $ Q_m$ for $1 \le k \le m$.       \\

Consider $ w^\E:=v^\E_{t_1}$ and note that $ w^\E \in
C^{1,\alpha}(\overline{\Omega_m})$ is a weak solution of
\begin{equation} \label{w_deri}
 \left\{
\begin{array}{rrl}
-\Delta_t w^\E  - c^\E(t) \cdot \nabla_t w^\E + \left( 1 + \frac{n_1-1}{(t_1+\E)^2}
\right) w^\E &=& h_{t_1} \qquad \hfill \mbox{ in } Q_m, \\
 w^\E &=& 0 \qquad \hfill  \mbox{ on } \partial_1 Q_m, \\
 w^\E_{t_i} &=&  0 \qquad \hfill  \mbox{ on } \partial_i Q_m, \; \; \mbox{ for $ 2
\le i \le m$}
\end{array}
\right.
\end{equation} where $ \partial_i Q_m:=\{ t \in \partial Q_m: t_i \in \{0,1\} \mbox{
and for } \forall j \neq i \mbox{ we have } t_j \notin \{0,1\} \}$.  Note there are
no issues regarding the boundary conditions since one has enough regularity to pass
the required derivatives on the boundary.  Note that the right hand side of
(\ref{w_deri}) is nonnegative.    Note a weak solution $w^\E$ of (\ref{w_deri})
satisfies

\begin{equation} \label{weak_ppp}
\int_{Q_m} \left( \nabla_t w^\E \cdot \nabla_t \eta + \left\{ 1 +
\frac{n_1-1}{(t_1+\E)^2 } \right\} w^\E \eta \right) d \mu_m^\E(t) = \int_{Q_m}
h_{t_1} \eta d \mu_m^\E(t), \forall \eta \in X,
\end{equation} where $h_{t_1} \ge 0$ in $Q_m$ and where $X:=\{ \eta \in H^1(Q_m):
\eta = 0 \mbox{ on } \partial_1 Q_m \}$.   In particular we can take $ \eta=(w^\E)^-
\in X$ (the negative part of $ w^\E$ to see that
\[ \int_{Q_m} \left(| \nabla_t (w^\E)^-|^2 + \left\{ 1 + \frac{n_1-1}{(t_1+\E)^2 }
\right\} ((w^\E)^-)^2\right) d \mu_m^\E(t) = - \int_{Q_m} h_{t_1} (w^\E)^- d
\mu_m^\E(t) \le 0 \] and hence $ (w^\E)^-=0$ a.e. in $ Q_m$.  After doing this for
each $ 1 \le k \le m$ we can conclude that $ v^\E \in K_m$.  \\

\noindent
 {\it Step 3.}  We now need to send $ \E \searrow 0$.  First recall that $v^\E$ satisfies
(\ref{weak_reg}) and whose weak formulation is given by

\begin{equation} \label{weak_ep}
\int_{Q_m} \left( \nabla_t v^\E\cdot \nabla_t  \eta +  v^\E \eta \right) d
\mu_m^\E(t) = \int_{Q_m} h  \eta d \mu_m^\E(t), \forall \eta \in H^1(Q_m).
\end{equation}

Taking $ \eta=v^\E$ one easily obtains a bound on $\{ v^\E \}_{0<\E \le 1}$ in $Y_m$
and hence after passing to a suitable subsequence one can assume that $ v^\E
\to\tilde{v}$ weakly in $Y_m$ and $ \tilde{v}$ satisfies (\ref{weak_u}) with
$v$ replaced with $ \tilde{v}$.   But recalling $v$ is already a solution of
(\ref{weak_u}) and since the solution is unique we see that $ v=\tilde{v} \in K_m$
after noting that since $K_m$ is convex and closed in $Y_m$ and hence it is weakly
closed.    \\

We now show that $v$ has added regularity.  To do this we obtains bounds on $ \{v^\E
\}_\E$ independent of $ 0<\E<1$.  Set $ \widehat{w}^\E(t):=\frac{w^\E(t)}{t_1+\E}=
\frac{v^\E_{t_1}(t)}{t_1+\E}$ and using (\ref{w_deri}) we see that $ \widehat{w}^\E
\in C^{1,\alpha}(\oo{Q_m})$ is a weak solution of

\begin{equation} \label{w_max}
 \left\{
\begin{array}{rrl}
-\Delta_t \widehat{w}^\E - b^\E(t) \cdot \nabla_t \widehat{w}^\E + \widehat{w}^\E
&=& \frac{h_{t_1} }{t_1+\E} \qquad \hfill \mbox{ in } Q_m, \\
 \widehat{w}^\E &=& 0 \qquad \hfill  \mbox{ on } \partial_1 Q_m, \\
 \widehat{w}^\E_{t_i} &=&  0 \qquad \hfill  \mbox{ on } \partial_i Q_m, \; \; \mbox{
for $ 2 \le i \le m$}
\end{array}
\right.
\end{equation}
 where $b^\E(t)=\big (b_1^\E(t),..., b_m^\E(t) \big)$ with  $b^\E_k(t):= \frac{n_k-1 + \alpha_k}{t_k+\E}$ where $ \alpha_1=2$
and $ \alpha_k=0$ for $ 2 \le k \le m$.   Using suitable reflections one can again
show that $ \widehat{w}^\E \in C^{2,\alpha} (\oo{Q_m})$ for each fixed $ 0<\E<1$.
Note that we have $ 0 \le \widehat{w}^\E$ and we now proceed to show that
$\widehat{w}^\E$ is bounded independently of $ 0<\E<1$.  So let
 $ t^0=(t^0_1,...,t^0_m) \in \oo{Q_m}$ be such that $ \widehat{w}^\E(t^0)=\sup_{Q_m} \widehat{w}^\E$
and we can assume $\widehat{w}^\E(t^0)>0$.
  If $ t^0 \in Q_m$ then we have $ \Delta_t w^\E(t^0) \le 0$. Now suppose that $
t^0 \in \partial Q_m$.   Note from the boundary condition we have $ t^0 \notin
\oo{ \partial_1 Q_m}$.  Set $ I:=\{ 1 \le i \le m: \exists \delta>0 \mbox{ such
that } t^0 + \tau e_i \in \overline{Q_m} \mbox{ for $ |\tau|<\delta $ } \}$ where
$e_i$ is the $i^{th}$ standard basis vector in $ \IR^m$ and let $J:=\{ 2,3,...,m\}
\backslash I$;  note that $ 1 $ is always an element of $I$.
  Then note that since $\widehat{w}^\E$ has a maximum at $t^0$ we
have $ \widehat{w}^\E_{t_i}(t^0)=0$ and $ \widehat{w}^\E_{t_i t_i}(t^0) \le 0$ for
all $ i \in I$;  note here we are just using single variable calculus at a
interior maximum point.   Now note that for $i \in J$ we have $t^0 \in
\oo{\partial_i Q_m}$ and hence we have $ \widehat{w}^\E_{t_i}(t^0)=0$ after
considering the boundary condition.  Using the fact that $ w^\E(t_0)$ has a
maximum at $t^0$ and since $ \widehat{w}^\E_{t_i}(t^0)=0$  we must have $
\widehat{w}^\E_{t_i t_i}(t^0) \le 0$.  From this we obtain that
  \begin{eqnarray*}
   \widehat{w}^\E(t^0) &  = & \Delta_t \widehat{w}^\E(t^0) + b^\E(t^0) \cdot
\nabla_t \widehat{w}^\E(t^0) + \frac{h_{t_1}(t^0) }{t_1^0+\E} \\
    & \le & \frac{h_{t_1}(t^0) }{t_1^0+\E}    \\
    &=& \frac{ \left( h_{t_1}(t^0) - h_{t_1}( \oo{t} ) \right) }{ | t^0 - \oo{t}|}
\frac{ |t^0 - \oo{t}|}{t_1^0+\E} \quad \mbox{ where $
\oo{t}=(0,t^0_2,...,t^0_m)$} \\
     & \le &   \frac{\| h \|_{C^{1,1}} |t_1^0|}{t_1^0+\E} \le \| h\|_{C^{1,1}},
     \end{eqnarray*}
     and so we have shown $ 0 \le \frac{v^\E_{t_k}(t)}{t_k+\E} \le \| h\|_{C^{1,1}}$
in $Q_m$ for $k=1$ and, by  using the same argument, it also holds for all  $ 1 \le k \le m$.  We now return to
(\ref{loppe}) which we recall was $-\Delta_t \oo{v^\E} + \oo{v^\E} = \oo{h} +
\oo{c^\E}(t) \cdot \nabla_t \oo{v^\E} $ in $  D_1.$    Note that $\oo{c^\E}(t)
\cdot \nabla_t \oo{v^\E}$ is bounded independently of $ \E$ after considering
the bound on $ \frac{v^\E_{t_k}}{t_k+\E}$.  We can then apply elliptic
regularity to see that $ \oo{v^\E}$ is bounded in $ C^{1,\alpha}(D_2)$
(independently of $ \E$) for $0< \alpha<1$.  We can now argue as before to show
that $ \oo{c^\E} \cdot \nabla_t \oo{v^\E}$ is bounded in $C^{0,\alpha}(D_2)$
independently of $ \E$ and hence we can apply elliptic regularity theory to see
that $ \oo{v^\E}$ is bounded in $C^{2,\alpha}(D_3)$ independently of $ \E$.
To obtain global regularity we need to perform the even extension across the
outer boundaries.  Note, as mentioned before,  the even extension of $h$ will
now only by Lipschitz,  but this is sufficient to show that $ v^\E$ is bounded
in $ C^{2,\alpha}(\oo{Q_m})$ independently of $ \E$ and this gives us the
desired result.

\hfill $\Box$

\end{document}